\documentclass[11pt]{amsart}
\allowdisplaybreaks[1]

\usepackage{amsmath}
\usepackage{amsxtra}
\usepackage{amssymb}
\usepackage{amscd}
\usepackage{epsfig}

 \newcommand{\A}{\mathbb A}

 \newcommand{\bee}{\begin{equation}}
 \newcommand{\eee}{\end{equation}}
 
 \newcommand{\Lb}{\mbox {\boldmath ${\Lambda}$}}
 \newcommand{\Gb}{\mbox {\boldmath ${\Gamma}$}}

 \newcommand{\Lbs}{\mbox{\scriptsize\boldmath ${\Lambda}$}}

 \newcommand{\Pb}{\mbox {\bf P}}
 \newcommand{\Pbs}{\mbox {\scriptsize{\bf P}}}

\newcommand{\be}{\begin{eqnarray}}
\newcommand{\ee}{\end{eqnarray}}
\newcommand{\supp}{\mbox{\rm supp}}

\newcommand{\dens}{\mbox{\rm dens}}

\newcommand{\freq}{\mbox{\rm freq}}
\newcommand{\Vol}{\mbox{\rm Vol}}

\newcommand{\eps}{{\mbox{$\epsilon$}}}
\newcommand{\e}{{\varepsilon}}

\newcommand{\R}{{\mathbb R}}

\newcommand{\Z}{{\mathbb Z}}
\newcommand{\C}{{\mathbb C}}

\newcommand{\Ak}{{\mathcal A}}

\newcommand{\Dk}{{\mathcal D}}

\newcommand{\Pk}{{\mathcal P}}

\newcommand{\Ok}{{\mathcal O}}
\newcommand{\Sk}{{\mathcal S}}
\newcommand{\Tk}{{\mathcal T}}

\newcommand{\dist}{\mbox{\rm dist}}

\newcommand{\Lam}{{\Lambda}}

\newcommand{\Gam}{\Gamma}
\newcommand{\om}{\omega}

 \newtheorem{theorem}{Theorem}[section]
 \newtheorem{lemma}[theorem]{Lemma}
 \newtheorem{prop}[theorem]{Proposition}
 \newtheorem{cor}[theorem]{Corollary}

 \newtheorem{defi}[theorem]{Definition}
 
 \newtheorem{remark}[theorem]{Remark}

\numberwithin{equation}{section}

\begin{document}

\title[Substitution Delone Sets with Pure Point Spectrum are Inter Model Sets]{Substitution Delone Sets with Pure Point Spectrum are Inter Model Sets}

\author{Jeong-Yup Lee}
\address{ KIAS 207-43 Cheongnyangni 2-dong, Dongdaemun-gu, Seoul 130-722, Korea
\hspace*{12pt}} \email{jylee@kias.re.kr or
jeongyuplee@yahoo.co.kr}

\date{\today}

\thanks{
\indent The author
acknowledges support from the NSERC post-doctoral fellowship and
KIAS research fellowship.}

\begin{abstract}
The paper establishes an equivalence between pure point
diffraction and certain types of model sets, called inter model
sets, in the context of substitution point sets and substitution
tilings. The key ingredients are a new type of coincidence
condition in substitution point sets, which we call algebraic
coincidence, and the use of a recent characterization of model
sets through dynamical systems associated with the point sets or
tilings.
\end{abstract}

\maketitle

\section{Introduction}

\noindent In the study of aperiodic order, there has been
considerable interest in understanding the structure of point sets
which have pure point diffraction spectrum. The pure point peaks
in the diffraction spectrum are indicative of a highly ordered
structure of the point sets, and the pure point diffraction
spectrum has provided a new viewpoint of looking at ordered
structures. Independent of the diffraction spectrum, there is the
spectral theory of dynamical system that arises as the completion
of the translational orbit of a point set in the standard
Radin-Wolff type topology. It turns out that pure point dynamical
and diffraction spectra are equivalent in quite a general setting
(see \cite{LMS1}, \cite{Gouere} and \cite{BL}). Especially, in
substitution point sets, the two types of pure point spectra are
equivalent, though not equal usually.

There is a large class of point sets which come from cutting and
projecting a lattice in a ``higher dimensional'' space $\R^d
\times H$ into the two lower dimension spaces $\R^d$ and $H$,
where $H$ is a locally compact Abelian group. Discrete point sets
in $\R^d$ are obtained by restricting the projection, which maps
from $\R^d \times H$ to $\R^d$, to some part of a lattice lying in
a cylinder of the form $\R^d \times W$, where the window $W
\subset H$ has non-empty interior and compact closure. If a
discrete point set in $\R^d$ comes from this projection and its
window has the boundary of measure zero, we call it a regular
model set (see Sec.\,\ref{prelim} for precise definitions).
This has provided a general way of obtaining point sets which have
the property of pure point dynamical spectrum (see \cite{BMS},
\cite{LM1} and \cite{LMS2}). The inversion problem, that is,
determining the structure of a discrete set knowing that it is
pure point diffractive, is in general impossible to solve.
However, with added ingredients it is possible to infer
information about the nature of a set from pure point
diffractivity. One such piece of information that we can now
determine, and that is the aim of this paper to prove, is the
structural type of the diffracting set - namely that it is a model
set, as long as we have the added ingredient of a primitive
substitution.

There is a recent characterization of model sets through the use
of dynamical systems associated with point sets in \cite{BLM} and
with multi-colour point sets in \cite{LM2}.
The notion of inter model sets is introduced in \cite{LM2} (under
the name of model sets) and \cite{BLM} as a model set satisfying a
topological condition which is less restrictive than the boundary
condition of a regular model set. In this paper we consider the
inter model sets (see Def.\,\ref{def-reg-model-set}) and show the
equivalence between inter model sets and pure point dynamical
spectrum in the context of primitive substitution point sets.

The main ingredient that establishes the connection between inter
model set and pure point dynamical spectrum is algebraic
coincidence. In the literature there are many types of
coincidences for substitution point sets and tilings which are
equivalent to the property of pure point dynamical spectrum of
these sets (see \cite{Dek}, \cite{BD}, \cite{soltil}, \cite{LMS2},
\cite{BK}). For example, in the class of constant length symbolic
substitutions Dekking's coincidence condition is well-known. It
says the following: suppose that $A = \{ a_1, \dots, a_m \} $ is a
finite alphabet with associated set of words $A^*$, and we are
given a primitive substitution $\sigma : A \to A^*$ for which the
length $l$ of each word $\sigma(a_i)$ is same, so called
``constant length substitution'', and the height is $1$. Then the
associated substitution dynamical system has pure point spectrum
if and only if it admits a coincidence, in the sense that there is
$k \in \Z_+$ such that $k \le l^n$ for some $n$ and the $k$-th
letter of each word $\sigma^n(a_i)$, $i \le m$, is same. Although
the various types of the coincidences are defined in slightly
different ways, they fulfill a similar property for pure point
dynamical spectrum. In the case of substitution point sets on
lattices it is established through {\em modular} coincidence that
regular model set is necessary and sufficient for pure point
dynamical spectrum and it is shown that the modular coincidence is
computable \cite{LMS2}. However there are many examples of
substitution point sets and tilings whose underlying structures
are not on lattices. For general substitution point sets we
introduce in this paper a new type of coincidence called {\em
algebraic coincidence}, whose name comes from the algebraic
structure of the point sets, which makes the construction of
locally compact Abelian groups and cut and project schemes
possible. The Fibonacci substitution sequence is a well-known
example of this kind.

Although we are primarily interested in and dealing with Delone
multi-colour sets, we need to introduce tilings along the way. It
is often advantageous to work with tilings in getting spacial
properties of point sets. However it is not necessarily true that
a substitution point set can be represented by a tiling in such a
way that every point of one type point set is represented by a
tile of the same type. In \cite{lawa}, Lagarias and Wang have
given a sufficient condition for a substitution point set to be
represented by a substitution tiling while maintaining its
iteration rules. They call it a self-replicating Delone set. In
\cite{LMS2} it is shown that a repetitive substitution point set
can be represented by a substitution tiling when a concept called
legality of clusters applies. We will talk more about this passage
from point sets to tilings in
Subsec.\,\ref{substi-pointSet-tilings}. We make use of this
connection in order to derive properties that we need in
substitution point sets from substitution tilings.

\medskip

\noindent
The main theorem in this paper states (see Th.\,\ref{main-theorem}) :

\medskip

\noindent {\em {\bf Theorem} Let $\Lb$ be a primitive substitution
Delone multi-colour set in $\R^d$ such that every $\Lb$-cluster is
legal and $\Lb$ has finite local complexity. Then the following
are equivalent:
\begin{itemize}
\item[(1)] $\Lb$ has pure point dynamical spectrum;
\item[(2)] $\Lb$ admits an algebraic coincidence;
\item[(3)] $\Lb$ is an inter model multi-colour set.
\end{itemize}
}

\medskip
\noindent
There are several recent results in the area that we use as new ingredients for the proof of the main theorem. We briefly explain them here.

We can construct two dynamical hulls generated by a point set,
using two different topologies. One is a local topology, which
defines the closeness of point sets by agreement on large regions
around the origin up to small shifts, and the other is a global
topology, called the autocorrelation topology, which defines the
closeness of point sets looking at how much they agree
density-wise up to small shifts. In \cite{BLM} conditions under
which there exists a continuous mapping between the two dynamical
hulls are given. In \cite{LM2}, it is shown that this mapping
derives inter model sets from model sets which are projected from
open windows.

A Delone set $\Lambda$ is called a Meyer set if $\Lambda -
\Lambda$ is uniformly discrete. The Meyer property is used
significantly in the main theorem to show the equivalence between
overlap coincidence on substitution tilings and algebraic
coincidence on substitution point sets, as well as in \cite{BLM}
and \cite{LM2} which call upon. However, in \cite{LS} it is shown
that any substitution point set with pure point dynamical spectrum
necessarily has the Meyer property, so we do not need to assume it
additionally in the main theorem.

We consider two topologies on the group $L$ generated by the
translation vectors of a substitution point set, thereby
constructing two locally compact Abelian groups which may be used
for defining CPSs. We call one $Q$-topology and the other
$P_{\eps}$-topology. In general, these two topologies on
substitution point sets are different. But under the assumption of
pure point spectrum, they are equivalent. Since model sets are
always associated with their CPSs, it is important to notice from
which CPS a model set arises. The equivalence of the two
topologies gives us the same CPS, and this allows us to freely use
the related results in \cite{BLM} and \cite{LM2}.

\medskip

The proof of the theorem is spread over several sections. The
structure of it is as follows: On substitution tilings it is known
that overlap coincidence is a necessary and sufficient condition
for pure point dynamical spectrum \cite{soltil}. We introduce
algebraic coincidence in substitution point sets, which is a
concept parallel to overlap coincidence, and in
Subsec.\,\ref{overlap-algebraic-coincidence} show the equivalence
between (1) and (2) of the theorem.
Prop.\,\ref{overlap-to-algebraic} plays an important role in
making connection between algebraic coincidence and pure point
spectrum. In Sec.\,\ref{construction-CPS}, we assume algebraic
coincidence in $\Lb$ and construct a CPS whose internal space is a
completion of a topological group $L$ with the $Q$-topology. We
then show that there exists a Delone multi-colour point set $\Gb$,
in a local dynamical hull generated by $\Lb$ which is a model set
with an open window in the CPS. We consider another topology
($P_{\eps}$-topology) on $L$ relative to which $L$ becomes a
topological group and show in Subsec.\,\ref{two-topology-on-L}
that the two topological spaces $L$ are in fact isomorphic. So
both topologies lead to the same completion of $L$. It is a
locally compact Abelian group and we can construct a CPS taking
this completed space as an internal space. In
Subsec.\,\ref{algebraic-coincidence-modelsets} we apply the
results of \cite{BLM} and \cite{LM2}, which are associated with
the $P_{\eps}$-topology, so that we get a condition for $\Lb$ to
be an inter model multi-colour set. We observe that algebraic
coincidence is sufficient for that condition to be fulfilled in
substitution point sets. We also show that the algebraic
coincidence is necessary to obtain the inter model multi-colour
set in Sec.\,\ref{InterModelsets-to-algCoincidence}.

The paper concludes with some unresolved questions (particularly
on the nature of the boundaries of the windows of the inter model
sets appearing in the theorem) and outlook for future work.

\medskip

\section{Preliminaries} \label{prelim}

\noindent Much of the terminology being introduced in this section
is standard and defined precisely in \cite{LMS2}. We refer the
reader to \cite{LMS2} for more detailed definitions and to
\cite{LP} for the standard concepts of discrete geometry in the
aperiodic setting.

\subsection{Delone multi-colour sets}

\noindent
A {\em multi-colour set} or {\em $m$-multi-colour set} in $\R^d$ is a
subset $\Lb = \Lam_1 \times \dots \times \Lam_m
\subset \R^d \times \dots \times \R^d$ \; ($m$ copies)
where $\Lam_i \subset \R^d$. We also write
$\Lb = (\Lam_1, \dots, \Lam_m) = (\Lam_i)_{i\le m}$.
Recall that a Delone set is a relatively dense and uniformly discrete subset of $\R^d$.
We say that $\Lb=(\Lambda_i)_{i\le m}$ is a {\em Delone multi-colour set} in $\R^d$ if
each $\Lambda_i$ is Delone and $\supp(\Lb):=\bigcup_{i=1}^m \Lambda_i \subset \R^d$ is Delone.
A {\em cluster} of $\Lb$ is, by definition,
a family $\Pb = (P_i)_{i\le m}$ where $P_i \subset \Lambda_i$ is
finite for all $i\le m$.
Many of the clusters that we consider have the form
$ \Lb \cap A := (\Lambda_i \cap A )_{i\le m}$, for a bounded set
 $A\subset \R^d$.
The translate of a cluster $\mbox{\bf P}$ by $x \in \R^d$ is
$x + \Pb = (x + P_i)_{i\le m}$.
We say that two clusters $\Pb$ and $\Pb'$ are {\em translationally equivalent}
if $\Pb=x+\Pb'$ for some $x \in \R^d$.
For any two Delone $m$-multi-colour sets $\Lb$ and $\Gb$,
we define $\Lb \, \cap \, \Gb = (\Lam_i \, \cap \, \Gam_i)_{i \le m}$ and
$\Lb \,\triangle \, \Gb = (\Lam_i \,  \triangle \, \Gam_i)_{i \le m}$, where $\Lam_i \, \triangle \, \Gam_i = (\Lam_i \, \backslash \, \Gam_i) \cup
(\Gam_i \, \backslash \, \Lam_i)$.
We write $B_R(y)$ for the open ball of radius $R$
centered at $y$ and use also $B_R$ for $B_R(0)$.
We define $\Xi(\Lb) := \bigcup_{i \le m}(\Lam_i - \Lam_i)$.
We say that $\Lam \subset \R^d$ is a {\em Meyer set} if it is a Delone set and $\Lam - \Lam \subset \Lam + F$ for some finite set $F$, equivalently, if it is a Delone set and $\Lam - \Lam$ is uniformly discrete (see \cite{Lag}, \cite{Moody1}).
We say $\Lb = (\Lam_i)_{i \le m}$ a {\em Meyer multi-colour set} if each component $\Lam_i$, $i \le m$, is a Meyer set and
$\supp(\Lb)$ is a Meyer set. A multi-colour set $\Lb$ is said to be {\em non-periodic} if there is no non-zero $x \in \R^d$ such that $\Lb + x = \Lb$.

A {\em cut and project scheme} (CPS) consists of a collection of spaces and mappings as follows;
\be  \label{CPS-definition}
\begin{array}{ccccc}
 \R^{d} & \stackrel{\pi_{1}}{\longleftarrow} & \R^{d} \times H & \stackrel{\pi_{2}}
{\longrightarrow} & H \\
 && \bigcup \\
 && \widetilde{L}
\end{array}
\ee
where $\R^{d}$ is a real Euclidean space, $H$ is some locally
compact Abelian group, $\pi_{1}$ and $\pi_{2}$ are the canonical projections,
$ \widetilde{L} \subset {\R^{d}
\times H}$ is a lattice,  i.e.\
a discrete subgroup for which the quotient group
$(\R^{d} \times H) / \widetilde{L}$ is
compact, $\pi_{1}|_{ \widetilde{L}}$ is injective,
and $\pi_{2}(\widetilde{L})$ is dense in $H$.

We call $\R^d$ a {\em physical space} and $H$ an {\em internal space}.
For a subset $V \subset H$, we denote $\Lambda(V) := \{\pi_1(x) \in \R^d : x \in \widetilde{L},
\pi_2(x) \in V \}$. We call the subset $V$ a {\em window} of $\Lambda(V)$.
A {\em model set} in $\R^d$ is a subset $\Gam$ of $\R^d$ for which $\Gam = \Lam(V)$ where
$V \subset H$ has non-empty interior and compact closure. The model set $\Gamma$ is
{\em regular} if the boundary
$\partial W = \overline{W} \backslash {W}^{\circ}$ of $W$
is of (Haar) measure $0$.

\begin{defi} \label{def-reg-model-set}
{\em An {\em inter model set} is a subset $\Gamma$ of $\R^{d}$ for which s + $\Lambda(W^{\circ}) \subset \Gamma \subset s+ \Lambda(W)$ for some $s \in \R^d$, where $W$ is compact in $H$ and $W = \overline{W^{\circ}} \neq \emptyset$,
with respect to CPS (\ref{CPS-definition}).}
\end{defi}

We say that $\Gb$ is a {\em model multi-colour set} (resp. {\em inter model multi-colour set}) if each $\Gamma_i$ is a model set (resp. inter model set) with respect to the same CPS (see \cite{Moody2}, \cite{LM2} and \cite{BLM} for more about model sets).
One should note here that since $\pi_2$ need not be $1-1$ on $\widetilde{L}$ in CPS, the notion of  inter model set, which is hemmed in between two such sets differing only by points on the boundary of the window $W$, arises naturally. When it is important to note which CPS a model set arises from, we will explicitly mention the CPS.

\medskip

Let $\Lb$ be a Delone multi-colour set.
We say that $\Lb$ has {\em finite local complexity
(FLC)} if for every $R > 0$ there exists a finite set $Y \subset \supp(\Lb) = \bigcup_{i=1}^m \Lam_i$ such that
for all $x \in \supp(\Lb)$, there exists $y \in Y$ for which $B_R(x) \cap \Lb = (B_R(y) \cap \Lb) + (x - y)$.
Also we say that $\Lb$ is {\em repetitive} if for every compact set
$K \subset \R^d$, $\{t \in \R^d : \Lb \cap K = (t + \Lb) \cap K\}$ is relatively dense.
For a cluster $\Pb$ and a bounded set $A\subset \R^d$, let
$$
L_{\Pbs}(A) := \sharp\{x\in \R^d:\ x+\Pb \subset A\cap \Lb\},
$$
where $\sharp$ means the cardinality.
A {\em van Hove sequence} for $\R^d$ is a sequence
$\mathcal{F}=\{F_n\}_{n \ge 1}$ of bounded measurable subsets of
$\R^d$ satisfying
\be \label{Hove}
\lim_{n\to\infty} \Vol((\partial F_n)^{+r})/\Vol(F_n) = 0,~
\mbox{for all}~ r>0,
\ee
where $(\partial{F_n})^{+r} := \{x \in \R^d:\,\dist(x, \partial{F_n}) \le r\}$.
We define
$$\dens(\Lb) :=
\lim_{n \to \infty} \frac{\sharp(\Lb \cap F_n)}{\Vol(F_n)},
$$
if the limit exists.
We say that $\Lb$ has {\em uniform cluster frequencies{\footnote{Here we define UCF with one fixed van Hove sequence. However it is implicit from \cite{LMS1} that UCF for the Delone multi-colour set with FLC does not depend on the choice of van Hove sequence.}}} (UCF) relative to $\{F_n\}_{n \ge 1}$ if for any cluster $\Pb$, there exists the limit
$$
\freq(\Pb,\Lb) = \lim_{n\to \infty} \frac{L_{\Pbs}(x+F_n)}{\Vol(F_n)},
$$
uniformly in $x \in \R^d$.

Let $X_{\Lbs}$ be the collection of all Delone multi-colour sets each of whose
clusters is a translate of a $\Lb$-cluster. We introduce a metric
on Delone multi-colour sets in a simple variation of the standard way:
 for Delone multi-colour sets $\Lb_1$, $\Lb_2 \in X_{\Lbs}$,
\be \label{metric-multi-colour sets}
d(\Lb_1,\Lb_2) := \min\{\tilde{d}(\Lb_1,\Lb_2), 2^{-1/2}\}\, ,
\ee
where
\be
\tilde{d}(\Lb_1,\Lb_2)
&=&\mbox{inf} \{ \e > 0 : \exists~ x,y \in B_{\e}(0), \nonumber \\ \nonumber
&  & ~~~~~~~~~~ B_{1/{\e}}(0) \cap (-x + \Lb_1) = B_{1/{\e}}(0)
\cap (-y + \Lb_2) \}\,.
\ee
For the proof that $d$ is a metric, see \cite{LMS1}.
Observe that $X_{\Lbs} = \overline{\{-h + \Lb : h \in \R^d \}}$ where the closure is taken in the topology induced by the metric $d$.
For more general topology defined by uniformity, see \cite{LM2}, \cite{martin} and \cite{BLM}.
We have a natural action of $\R^d$ on the dynamical hull $X_{\Lbs}$ of $\Lb$ by translations which makes it a topological dynamical system $(X_{\Lbs},\R^d)$. With FLC, $X_{\Lbs}$ is a compact space.

\medskip

Let $(X_{\Lbs}, \mu, \R^d)$ be a measure preserving dynamical system.
We consider the associated group of unitary operators $\{T_x\}_{x\in \R^d}$ on
$L^2(X_{\Lbs},\mu)$:
$$
T_x g(\Lb') = g(-x + \Lb').
$$
Every $g \in L^2(X_{\Lbs},\mu)$ defines a function on $\R^d$ by
$x \mapsto \langle T_x g,g \rangle$.
This function is positive definite on $\R^d$, so its
Fourier transform is a positive measure $\sigma_g$ on $\R^d$ called the
{\em spectral measure} corresponding to $g$.
The dynamical system $(X_{\Lbs}, \mu, \R^d)$ is said to have {\em
pure point spectrum} if $\sigma_g$ is pure point for every $g \in
L^2(X_{\Lbs}, \mu)$. \footnote{We also say that $\Lb$ has pure
point spectrum if $\sigma_g$ is pure point for every $g \in
L^2(X_{\Lbs},\mu)$}. We recall that $g \in L^2(X_{\Lbs},\mu)$ is
an eigenfunction for the $\R^d$-action if for some $\alpha
=(\alpha_1,\ldots,\alpha_d) \in \R^d$,
$$
T_x g = e^{2 \pi i x\cdot \alpha} g
\ \ \ \mbox{for all}\ \ x\in \R^d,
$$
where $x\cdot \alpha$ is the standard inner product on $\R^d$.

\subsection{Tilings}

We begin with a set of types (or colours) $\{1,\ldots,m\}$, which
we fix once and for all. A {\em tile} in $\R^d$ is defined as a
pair $T=(A,i)$ where $A=\supp(T)$ (the support of $T$) is a
compact set in $\R^d$, which is the closure of its interior, and
$i=l(T)\in \{1,\ldots,m\}$ is the type of $T$. We let $g+T =
(g+A,i)$ for $g\in \R^d$. We say that a set $P$ of tiles is a {\em
patch} if the number of tiles in $P$ is finite and the tiles of
$P$ have mutually disjoint interiors. The {\em support of a patch}
is the union of the supports of the tiles that are in it. The {\em
translate of a patch} $P$ by $g\in \R^d$ is $g+P := \{g+T:\ T\in
P\}$. We say that two patches $P_1$ and $P_2$ are {\em
translationally equivalent} if $P_2 = g+P_1$ for some $g\in \R^d$.
A tiling of $\R^d$ is a set $\Tk$ of tiles such that $\R^d =
\bigcup \{\supp(T) : T \in \Tk\}$ and distinct tiles have disjoint
interiors. Given a tiling $\Tk$, a finite set of tiles of $\Tk$ is
called $\Tk$-patch. We define FLC, repetitivity, and uniform patch
frequencies (UPF), which is the analog of UCF, on tilings in the
same way as the corresponding properties on Delone multi-colour
sets. The types (or colours) of tiles on tilings have the same
concept as the colours of points on Delone multi-colour sets. We
always assume that any two $\Tk$-tiles with the same colour are
translationally equivalent (hence there are finitely many
$\Tk$-tiles up to translations).

For a subset $\mathcal{S}$ of a tiling and $A \subset \R^d$, we define
\[ \mathcal{S} \cap A := \{T \in \mathcal{S} : (\supp{(T)})^{\circ} \cap A \neq \emptyset\}
\]
and for tilings $\Tk$ and $\Tk'$, we use $\Tk \cap \Tk' \cap A$ for $(\Tk \cap \Tk') \cap A$.
For any patch $P$, we write $\Vol(P)$ for $\Vol( \bigcup\{\supp(T) : T \in P\})$.
Just as for point sets, for any $D \subset \Tk$, we define
\[
\dens(D) : =  \lim_{n \to \infty}\frac{\Vol(D \cap F_n)}{\Vol(F_n)} \ \  \mbox{and} \ \
 \freq(P, \Tk) :=\lim_{n\to \infty} \frac{L_{P}(F_n)}{\Vol(F_n)},
\]
if the limits exist.
Let $X_{\Tk}$ be the collection of all tilings each of whose patches is a translate of $\Tk$-patch.
We define a metric $d$ on tilings, given analogously to (\ref{metric-multi-colour sets}) for Delone multi-colour sets:
for tilings $\Tk$, $\Sk \in X_{\Tk}$,
\be \label{metric-tilings}
d(\Tk,\Sk) := \min\{\tilde{d}(\Tk,\Sk), 2^{-1/2}\}\, ,
\ee
where
\be
\tilde{d}(\Tk,\Sk)
&=&\mbox{inf} \{ \e > 0 : \exists~ x,y \in B_{\e}(0), \nonumber \\ \nonumber
&  & ~~~~~~~~~~  (-x + \Tk) \cap B_{1/{\e}}(0) = (-y + \Sk) \cap B_{1/{\e}}(0) \}\,.
\ee
We define the dynamical hull $(X_{\Tk}, \R^d)$ of $\Tk$ in the same way as of Delone multi-colour sets (see \cite{soltil}).
Also we have the equivalent notion of pure point spectrum on tilings.

\subsection{Substitutions} \label{substi-pointSet-tilings}

\noindent
\subsubsection{Substitutions on Delone multi-colour sets}

We say that a linear map $Q : \R^d \rightarrow \R^d$ is {\em expansive}
if there is
a $c > 1$ with
\be
e(Qx,Qy) \geq c \cdot e(x,y) \label{def-expan}
\ee
for all $x,y \in \R^d$
and some metric $e$ on $\R^d$ compatible with the standard topology.
This is equivalent to saying that all the eigenvalues of $Q$
 lie outside the closed unit disk in $\C$.

\begin{defi} \label{def-subst-mul}
{\em $\Lb = (\Lam_i)_{i\le m}$ is called a {\em
substitution Delone multi-colour set} if $\Lb$ is a Delone multi-colour set and
there exist an expansive map
$Q:\, \R^d\to \R^d$ and finite sets $\Dk_{ij}$ for $i,j\le m$ such that
\be \label{eq-sub}
\Lambda_i = \bigcup_{j=1}^m (Q \Lambda_j + \Dk_{ij}),\ \ \ i \le m,
\ee
where the unions on the right-hand side are disjoint.}
\end{defi}

We say that the
substitution Delone multi-colour set is {\em primitive} if the corresponding
substitution matrix $S$, with $S_{ij}= \sharp (\Dk_{ij})$, is primitive.
For any given substitution Delone multi-colour set $\Lb = (\Lambda_i)_{i \le m}$,
we define $\Phi_{ij} = \{ f : x \mapsto Qx + a \, : \,a \in \Dk_{ij}\}$.
Then $\Phi_{ij}(\Lam_j) = Q \Lam_j + \Dk_{ij}$, where $i \le m$. We define $\Phi$ an $m \times m$ array for which each entry is $\Phi_{ij}$, and call $\Phi$ a {\em matrix function system (MFS)} for the substitution.
For any $k \in \Z_+$ and $x \in \Lam_j$ with $j \le m$, we let
$\Phi^k (x) =  \Phi^{k-1}((\Phi_{ij}(x))_{i \le m})$.
For any $k \in \Z_+$, $\Phi^k (\Lb) = \Lb$ and
$\Phi^k (\Lam_j) = \bigcup_{i \le m}(Q^k \Lam_j + (\Dk^k)_{ij})$ where
\[(\Dk^k)_{ij} = \bigcup_{n_1,n_2,\dots,n_{(k-1)} \le m}
(\Dk_{in_1} + Q \Dk_{n_1 n_2} + \cdots + Q^{k-1} \Dk_{n_{(k-1)} j}).
\]
We say that a cluster $\Pb$ is {\em legal} if it is a translate of a subcluster of
a cluster generated from one point of $\Lb$, i.e. $a + \Pb \subset \Phi^k (x)$ for some $k \in \Z_+$, $a \in \R^d$ and $x \in \Lb$.

\subsubsection{Substitutions on tilings}

\begin{defi}\label{def-subst}
{\em Let $\Ak = \{T_1,\ldots,T_m\}$ be a finite set of tiles in $\R^d$
such that $T_i=(A_i,i)$; we will call them {\em prototiles}.
Denote by $\Pk_{\Ak}$ the set of
patches made of tiles each of which is a translate of one of $T_i$'s.
We say that $\omega: \Ak \to \Pk_{\Ak}$ is a {\em tile-substitution} (or simply
{\em substitution}) with
expansive map $Q$ if there exist finite sets $\Dk_{ij}\subset \R^d$ for
$i,j \le m$, such that
\begin{equation}
\om(T_j)=
\{u+T_i:\ u\in \Dk_{ij},\ i=1,\ldots,m\}
\label{subdiv}
\end{equation}
with
$$
Q A_j = \bigcup_{i=1}^m (\Dk_{ij}+A_i) \ \ \  \mbox{for} \  j\le m.
$$
Here all sets in the right-hand side must have disjoint interiors;
it is possible for some of the $\Dk_{ij}$ to be empty.}
\end{defi}

Note that $Q A_j = \supp (\om(T_j)) = Q \supp(T_j)$.
The substitution (\ref{subdiv}) is extended to all translates of prototiles by
\be \label{tile-equivalence}
\om(x+T_j)= Q x + \om(T_j),
\ee
in particular,
\be \label{supp-of-iterated-tile}
\supp (\om(x + T_j)) = \supp(Qx + \om(T_j)) = Qx + Q \supp(T_j) = Q(x + \supp(T_j)),
\ee
and to patches and tilings by
$\om(P)=\cup\{\om(T):\ T\in P\}$.
The substitution $\om$ can be iterated, producing larger and larger patches
$\om^k(T_j)$.

We define the substitution matrix and primitivity of $\omega$ in the similar way as in substitution Delone multi-colour sets. We say that $\Tk$ is a {\em substitution tiling} if $\Tk$ is a tiling and $\om(\Tk) = \Tk$ with some substitution $\om$.
We say that a patch $P$ is {\em legal} if it is a translate
of a subpatch of $\om^k(T_i)$ for some $i\le m$ and $k\ge 1$. This is the analog of a legal cluster on Delone multi-colour sets.

\subsubsection{Representability of $\Lb$ as a tiling} \label{RepresentabilityAsTiling}

Let $\Lb$ be a substitution Delone multi-colour set. One can set
up an {\em adjoint system} of equations \be \label{eq-til} Q A_j =
\bigcup_{i=1}^m (\Dk_{ij} + A_i),\ \ \ j \le m \ee from the
equation (\ref{eq-sub}). It is known that (\ref{eq-til}) always
has a unique solution for which $\Ak = \{A_1, \dots, A_m\}$ is a
family of non-empty compact sets of $\R^d$ (see for example
\cite{BM1}, Prop.\,1.3). It is proved in \cite[Th.\,2.4 and
Th.\,5.5]{lawa} that if all the sets $A_i$ from (\ref{eq-til})
have non-empty interiors and, moreover, each $A_i$ is the closure
of its interior. We say that $\Lb$ is {\em representable} (by
tiles) if $\Lb + \Ak := \{x + T_i :\ x\in \Lambda_i,\ i \le m\}$
is a tiling of $\R^d$, where $T_i = (A_i,i)$, $i \le m$, and
$A_i$'s arise from the solution to the adjoint system
(\ref{eq-til}) and $\Ak = \{T_i : i\le m\} $. One can define a
tile-substitution $\omega$ satisfying $\omega(\Lb + \Ak) = \Lb +
\Ak$ from (\ref{eq-til}). So $\Lb + \Ak$ is a substitution tiling.
We call $\Lb + \Ak$ the {\em associated substitution tiling} of
$\Lb$. In \cite[Lemma 3.2]{lawa} it is shown that if $\Lb$ is a
substitution Delone multi-colour set, then there is a finite
multi-colour set (cluster) $\Pb \subset \Lb$ for which
$\Phi^{n-1}(\Pb) \subset \Phi^n(\Pb)$ for $n \ge 1$ and $\Lb =
\lim_{n \to \infty} \Phi^n (\Pb)$. We call such a multi-colour set
$\Pb$ a {\em generating set} for $\Lb$.

\begin{theorem}\cite{LMS2} \label{legal-rep}
Let $\Lb$ be a repetitive primitive substitution Delone multi-colour set. Then every
$\Lb$-cluster is legal if and only if $\Lb$ is representable.
\end{theorem}

Note that in order to check that every $\Lb$-cluster is legal, it suffices to check if
a cluster containing a finite generating set for $\Lb$ is legal (see \cite{LMS2}).

\begin{remark}
{\em Throughout this paper we are mainly interested in primitive
substitution Delone multi-colour sets $\Lb$ such that every
$\Lb$-cluster is legal. Since $\Lb$ is representable, we will
often identify $\Lb$ with the associated substitution tilings $\Tk
= \Lb + \Ak = \{ x + T_i : x \in \Lam_i, i \le m\}$, where $T_i$'s
are the tiles arising from the solution of the adjoint system of
equations.}
\end{remark}

For a primitive representable substitution Delone multi-colour set $\Lb$, the dynamical system $(X_{\Lbs}, \R^d)$ is unique ergodic. Similarly the dynamical system $(X_{\Tk}, \R^d)$ of any primitive substitution tiling $\Tk$ is uniquely ergodic \cite{LMS2}.

\section{Algebraic coincidence and pure point spectrum} \label{overlap-algebraic-coincidence}

\noindent
For substitution tilings overlap coincidence was introduced in \cite{soltil}. The overlap coincidence has been a central concept connecting the asymptotic behavior of $Q$-iterates of almost-periods of tilings with pure pointedness of the spectrum of dynamical systems generated by the tilings (see \cite{soltil} and \cite{LMS2}).
Here we introduce a new coincidence condition for substitution Delone multi-colour sets and show that it is equivalent to the overlap coincidence when the substitution Delone multi-colour sets are representable.

\medskip
Let $\Tk$ be a tiling and $\Xi(\Tk)$ be the set of translation
vectors between $\Tk$-tiles of the same type: \be \label{def-xi}
\Xi(\Tk) := \{x\in \R^d:\ \exists \,T,T' \in \Tk, \ T'=x+T \}. \ee
Since $\Tk$ has the inflation symmetry with the expansive map $Q$,
we have that $Q\Xi(\Tk) \subset \Xi(\Tk)$. Note also that
$\Xi(\Tk) = -\Xi(\Tk)$. If $\Tk=\Lb + \Ak$ is an associated
substitution tiling of $\Lb$, then $\Xi(\Tk) = \bigcup_{i=1}^m
(\Lambda_i-\Lambda_i)$. Recall that $\omega{\Tk} = \Tk$, where
$\omega$ is the tile-substitution coming from the adjoint system
of equations.

\begin{defi}\label{def-overlap} {\em Let $\Tk$ be a tiling. A triple $(T,y,S)$,
with $T,S \in \Tk$ and $y \in \Xi(\Tk)$, is called an {\em overlap}
if $\supp(y+T) \, \cap \, \supp(S)$ has non-empty interior.
We say that two overlaps $(T,y,S)$ and $(T',y',S')$ are {\em equivalent} if
for some $g\in \R^d$ we have $y+T = g+y'+T',\ S = g+S'$. Denote by
$[(T,y,S)]$ the equivalence class of an overlap. An overlap $(T,y,S)$ is a
{\em coincidence} if $y+T = S$. The {\em support} of an overlap
$(T,y,S)$ is $\supp(T,y,S) = \supp(y+T) \,\cap \, \supp(S)$.}
\end{defi}

Let $\Ok=(T,y,S)$ be an overlap. Recall that for
a tile-substitution $\om$,
$\om(y+T) = Qy+\om(T)$ is a patch of $Qy+\Tk$, and $\om(S)$ is a $\Tk$-patch,
and moreover,
$$
\supp(Qy+\om(T)) \cap \supp(\om(S)) = Q(\supp(T,y,S)).
$$
For each $l \in \Z_+$,
$$Q^l (\mathcal{O}) = \{(T',Q^l y,S') : T'\in \om^l(T), S'\in \om^l(S),
\supp(Q^l y+T') \, \cap \, \supp(S') \neq \emptyset \}.$$

\begin{defi}
{\em We say that a substitution tiling $\Tk$ admits an {\em overlap coincidence} if there exists $l \in \Z_+$ such that for each overlap $\mathcal{O}$ in $\Tk$, $Q^l (\mathcal{O})$ contains a coincidence.}
\end{defi}

\begin{theorem} \cite[Th.\,4.7 and Lemma\,A.9]{LMS2} \label{Ppd-overlap-thm}
Let $\Tk$ be a repetitive fixed point of a primitive substitution such that
$\Xi(\Tk)$ is a Meyer set.
Then $(X_{\Tk}, \R^d, \mu)$ has a pure point dynamical spectrum if and only if $\Tk$ admits an overlap coincidence.
\end{theorem}

\begin{theorem} \cite{LS} \label{purePoint-Meyer}
Let $\Lb$ be a primitive substitution Delone multi-colour set such that every $\Lb$-cluster is legal and $\Lb$ has FLC.
Suppose that $(X_{\Lbs}, \R^d, \mu)$ has a pure point dynamical spectrum. Then
$\Lambda = \bigcup_{i \le m} \Lambda_i$ and $\Xi(\Lb)$ are Meyer sets.
\end{theorem}

From the result of Th.\,\ref{purePoint-Meyer}, we can get the following corollary of Th.\,\ref{Ppd-overlap-thm} dropping the Meyer condition.

\begin{cor} \label{purePoint-OverlapCoincidence}
Let $\Tk$ be a repetitive fixed point of a primitive substitution with FLC. Then $(X_{\Tk}, \R^d, \mu)$ has a pure point dynamical spectrum if and only if $\Tk$ admits an overlap coincidence.
\end{cor}

\noindent {\sc Proof.} We first note that when we replace the third condition of \cite[Lemma\,A.9]{LMS2} by the overlap coincidence that we define here, Lemma\,A.9 in \cite{LMS2} holds without the assumption of the Meyer property. So applying \cite[Th.\,6.1]{soltil}, we can see that the necessity direction holds. The sufficiency follows from Th.\,\ref{Ppd-overlap-thm} and Th.\,\ref{purePoint-Meyer}.   \qed

\medskip
Applying Th.\,\ref{purePoint-Meyer} to Cor.\,\ref{purePoint-OverlapCoincidence}, we observe the following.

\begin{cor} \label{overlapCoincidence-MeyerSet}
Let $\Tk$ be a repetitive fixed point of a primitive substitution with FLC. If $\Tk$ admits an overlap coincidence, then $\Xi(\Tk)$ is a Meyer set. \qed
\end{cor}

\begin{lemma} \cite[Lemma\,A.8]{LMS2}\label{lem-finite}
Let $\Tk$ be a tiling
such that $\Xi(\Tk)$ is a Meyer set. Then the number of equivalence classes
of overlaps for $\Tk$ is finite.
\end{lemma}

\begin{defi}
{\em Let $\Lb$ be a primitive substitution Delone multi-colour set with expansive map $Q$. We say that $\Lb$ admits an {\em algebraic coincidence} if there exist $M \in \Z_+$ and $\xi \in \Lam_i$ for some $i \le m$ such that $\xi + Q^M \Xi(\Lb) \subset \Lam_i$. We say that $\Lb$ admits an {\em algebraic coincidence at $\xi$}, when we need to emphasize the role of $\xi$ for the algebraic coincidence.}
\end{defi}

\begin{lemma} \label{alg-coincidence-Meyer}
Let $\Lb$ be a primitive substitution Delone multi-colour set with expansive map $Q$.
If $\Lb$ admits an algebraic coincidence, then $\Xi(\Lb)$ is a Meyer set and thus $\Lb$ has FLC.
\end{lemma}

\noindent {\sc Proof.} By the assumption, there exist $M \in \Z_+$
and $\xi \in \Lam_i$ for some $i \le m$ such that $\xi + Q^M
\Xi(\Lb) \subset \Lam_i$. Since $\Lb$ is a Delone multi-colour
set, $Q^M \Xi(\Lb)$ is uniformly discrete. So $\Xi(\Lb)$ is
uniformly discrete. Note that $Q^M \Xi(\Lb) - Q^M \Xi(\Lb) \subset
\Lam_i - \Lam_i \subset \Xi(\Lb)$. Thus $\Xi(\Lb) - \Xi(\Lb)$ is
uniformly discrete, i.e. $\Xi(\Lb)$ is a Meyer set. Then it is
easy to see that $\Lb$ has FLC.  \qed

\medskip

We now connect overlap coincidence with algebraic coincidence.
For this connection, the condition of $\Xi(\Lb)$ being a Meyer set is strongly used.

\begin{prop} \label{overlap-to-algebraic}
Let $\Lb$ be a primitive substitution Delone multi-colour set such that every $\Lb$-cluster is legal
and $\Lb$ has FLC.
Suppose that the associated substitution tiling $\Tk = \Lb + \Ak$
admits an overlap coincidence. Then $\Lb$ admits an algebraic coincidence.
\end{prop}

\noindent {\sc Proof.} We first sketch the idea of the proof.
First, we note that $\Xi(\Tk)$ is a Meyer set by Cor\,\ref{overlapCoincidence-MeyerSet}.
From Lemma \ref{lem-finite}, there are only finitely many possible
overlaps in $\Tk$. We start with any nonempty patch $P$ in $\Tk$.
As $\Tk$ is translated by the vectors in $\Xi(\Tk)$ and intersects
with $\Tk$, we get certain configurations of overlaps on
$\supp(P)$. There are only finitely many possible configurations
of overlaps on $\supp(P)$ that arise this way. When each
configuration of overlaps on $\supp(P)$ is enlarged enough by
applying the substitution repeatedly, there are many coincidences
of overlaps occuring in the enlarged configuration. Ultimately
they cover most of the volume of the enlarged configuration, due
to the assumption of overlap coincidence (see \cite[Lemma\
A.9]{LMS2}). But the number of configurations of overlaps on
$\supp(P)$ stays same as we enlarge them. So when we intersect all
the coincidences of overlaps in the enlarged configuration for all
the translational vectors in $\Xi(\Tk)$, we get a {\em non-empty}
set. This implies that there exists at least one tile in $\Tk$
whose translations by the translational vectors of $Q^M \Xi(\Tk)$
are all in $\Tk$. This implies algebraic coincidence for $\Lb$.

Now we give a detailed proof.
We consider the associated substitution tiling $\Tk = \Lb + \Ak$ of $\Lb$ and choose any nonempty patch $P$ in $\Tk$. Note that $\Xi(\Tk) = \Xi(\Lb)$. Consider the collection of patches of translates of
$\Tk$ on $\supp(P)$.
$$\mathcal{H} = \{(\alpha +\Tk) \cap \supp(P) : \alpha \in \Xi(\Tk) \}.$$
Since $\Xi(\Tk)$ is a Meyer set, the number of translationally equivalent classes of overlaps for $\Tk$ is finite by Lemma \ref{lem-finite}.
It is important to note as a result of this that $\mathcal{H}$ consists of only finitely many patches.
Notice that this is more than just saying that there are finitely many translational classes of patches, which simply means the
FLC of $\Lb$.
Thus we can find
$\alpha_1, \dots, \alpha_K \in \Xi(\Tk)$ such that for any $\alpha \in \Xi(\Tk)$,
\be \label{finite-patch}
(\alpha + \Tk) \cap \supp(P) = (\alpha_k + \Tk) \cap \supp(P)  \ \ \  \mbox{for some} \ k \le K.
\ee
For any $n \in \Z_+$, we get
\be \label{infl-finite-patch}
\omega^n ((\alpha + \Tk) \cap \supp(P)) =  \omega^n ((\alpha_k + \Tk) \cap \supp(P)).
\ee
Looking at the set of (\ref{infl-finite-patch}) on the compact set $Q^n \supp(P)$, we get from $\omega^n \Tk = \Tk$ and
(\ref{tile-equivalence}) that
\[ (Q^n \alpha + \Tk) \cap Q^n \supp(P)  = (Q^n \alpha_k + \Tk) \cap Q^n \supp(P).\]
Then
\be \label{same-Q^nB-cluster}
\Tk \cap (Q^n \alpha + \Tk) \cap Q^n \supp(P)  = \Tk \cap (Q^n \alpha_k + \Tk)
\cap Q^n \supp(P).
\ee
Let
\be
D_{Q^{n} \alpha} :=  \Tk \cap (Q^{n} \alpha + \Tk) \cap Q^{n} \supp(P),
\ee
for any $\alpha \in \Xi(\Tk)$.
Then from (\ref{finite-patch}),
\be \label{D-setAll-D-setFinite}
\bigcap_{\alpha \in \Xi(\Tk)} D_{Q^{M} \alpha} = \bigcap_{k=1}^K D_{Q^{M} \alpha_k}
\ee
We claim that there exists $M \in \Z_+$ such that
\be \label{D-setFiniteNonempty}
 \bigcap_{k = 1}^K D_{Q^{M} \alpha_k} \neq \emptyset.
\ee
Then it implies that there exists $T \in \bigcap_{\alpha \in \Xi(\Tk)} D_{Q^{M} \alpha}$ from (\ref{D-setAll-D-setFinite})
and so for any $\alpha \in \Xi(\Tk)$ there exists $T' \in \Tk$ such that $Q^M \alpha + T' = T$.
Note that $ T = \xi + T_i $ and $T' = \eta + T_i $ \ where $\xi, \eta \in \Lam_i$ for some $i \le m$.
So $\xi - Q^M \alpha \in \Lam_i$ and
we can conclude that $\xi - Q^M \Xi(\Lb) = \xi + Q^M \Xi(\Lb) \subset \Lam_i$.
Therefore $\Lb$ admits an algebraic coincidence.

Now we give the proof of the claim (\ref{D-setFiniteNonempty}). We use a type of argument, leading to (\ref{calculation-non-overlaps}), that has been used before in \cite{soltil} and \cite{LMS2}.
However it is not in the form that we can make direct use here, so we discuss it again in the form that we need.
Let $V_0$ and $V_1$ be the minimal and maximal volumes of $\Tk$-tiles respectively.
Notice that $\supp(D_{Q^{n} \alpha})$ is the union of supports of coincidences in $Q^{n} \supp(P)$. It is easy to see that coincidence in $D_{Q^{n} \alpha}$ leads to other coincidences in
$D_{Q^{n+1} \alpha}$. Thus
$$ Q \supp(D_{Q^{n} \alpha}) \subset \supp(D_{Q^{n+1} \alpha}).$$
Since $\Tk$ admits an overlap coincidence, there exists $l \in \Z_+$ such that for each overlap $\mathcal{O}$ in $\Tk$, $Q^l (\mathcal{O})$ contains a coincidence.
Note that $\supp(D_{Q^{n+l} \alpha})$ has more support than $Q^l \supp(D_{Q^{n} \alpha})$ from the new coincidence occuring after $l$-step iterations of each non-coincidence overlap, and the volume of the support of the new coincidence is at least $V_0$ for
every $l$-step iteration of each non-coincidence overlap.
So we can get the following formula
\begin{eqnarray} \label{coincidence-volume-increase}
\lefteqn{\Vol(D_{Q^{n+l} \alpha} ) - |\det Q|^l \Vol(D_{Q^{n} \alpha})} \nonumber \\
&\ge & \frac{V_0}{{V_1} |\det Q|^l} \left( \Vol(Q^{n+l} \supp(P) ) - |\det Q|^l \Vol(D_{Q^{n} \alpha} ) \right).
\end{eqnarray}
Letting
$$b = 1 - \frac{V_0}{V_1 |\det Q|^l}$$
and using
$$\Vol(Q^{n+l} \supp(P)) = |\det Q|^l \Vol(Q^n \supp(P)) \,,$$
the inequality (\ref{coincidence-volume-increase}) becomes
\begin{eqnarray} \label{overlap-increase-formula}
1 - \frac{\Vol(D_{Q^{n+l} \alpha} )}{\Vol(Q^{n + l} \supp(P) )}
 \le  b \left( 1 - \frac{\Vol(D_{Q^{n} \alpha} )}{\Vol(Q^{n} \supp(P) )} \right).
\end{eqnarray}
For all $n = tl + s \ge 0$ where $t \in \Z_+$ and $0 \le s < l$, we obtain from (\ref{overlap-increase-formula})
\begin{eqnarray} \label{calculation-non-overlaps}
1 - \frac{\Vol(D_{Q^{n} \alpha})}{\Vol(Q^{n} \supp(P) )}
& \le & b^t \left( 1 - \frac{\Vol(D_{Q^{s} \alpha})}{\Vol(Q^{s} \supp(P) )} \right) \nonumber \\
& = & (b^{1/l})^{tl + s} \frac{1}{b^{s/l}}
\left( 1 - \frac{\Vol(D_{Q^{s} \alpha})}{\Vol(Q^{s} \supp(P) )} \right) \nonumber \\
& \le & r^n c \ \ \ \mbox{for some} \ r \in (0,1) \ \mbox{and}~ c > 0.
\end{eqnarray}
Thus for any $\epsilon > 0$, we can find $M \in \Z_+$ such that for any $1 \le k \le K$,
\[
1 - \frac{\Vol(D_{Q^{M} \alpha_k} )}{\Vol(Q^{M} \supp(P) )} < \epsilon .
\]
This implies that
\[ 1 - \frac{\Vol(\bigcap_{k=1}^K D_{Q^{M} \alpha_k} )}{\Vol(Q^M \supp(P) )}
< \epsilon K .
\]
Therefore for small $\epsilon > 0$
\be
\bigcap_{k=1}^K D_{Q^{M} \alpha_k} \neq \emptyset \,,   \nonumber
\ee
as we claimed in (\ref{D-setFiniteNonempty}).
\qed

\medskip

We show now the converse direction of Prop.\,\ref{overlap-to-algebraic}.
From Lemma\,\ref{alg-coincidence-Meyer}, we do not need to additionally assume the Meyer property of $\Xi(\Lb)$ in the following proposition.

\begin{prop} \label{algebraic-overlap-coincidence}
Let $\Lb$ be a primitive substitution Delone multi-colour set such that every $\Lb$-cluster is legal.
Suppose that $\Lb$ admits an algebraic coincidence. Then the associated substitution tiling
$\Tk =\Lb + \Ak$ admits an overlap coincidence.
\end{prop}

\noindent
{\sc Proof.}
Suppose that there exist $M \in \Z_+$ and $\xi \in \Lam_i$  such that
$\xi + Q^M \Xi(\Lb) \subset \Lam_i$ for some $i \le m$. Then
\be \label{Xi-closed-under-addition}
Q^{M} \Xi(\Lb) + Q^M \Xi(\Lb) & = & Q^M \Xi(\Lb) - Q^M \Xi(\Lb) \nonumber \\ \nonumber
& \subset & (\Lam_i - \xi) - (\Lam_i - \xi) \\
& = & \Lam_i - \Lam_i \subset \Xi(\Lb).
\ee
Thus $\xi + Q^{2M}\Xi(\Lb) + Q^{2M}\Xi(\Lb) \subset \Lam_i$.
For any overlap $\mathcal{O}=(R,y,S)$ with $R, S \in \Tk$, $y \in \Xi(\Lb)$, we can find
$M' = M'(\mathcal{O}) \ge 2M$ so that
\[\Tk \cap (Q^{M'}((y + \supp(R)) \cap \supp(S)) - Q^{M'} y)\]
contains at least one $\xi + Q^{2M} z + T_i$ with some $z \in \Lam_i - \Lam_i$, since $Q^{2M}(\Lam_i - \Lam_i)$ is relatively dense. Then
\[ \xi + Q^{2M} z + T_i \in \omega^{M'} (R) \,.\]
Note that $Q\Xi(\Lb) \subset \Xi(\Lb)$.
So
\[ \xi + Q^{2M}z + Q^{M'} y \subset \xi + Q^{2M}\Xi(\Lb) + Q^{2M}\Xi(\Lb) \subset \Lam_i \]
and
\[ \xi + Q^{2M} z + Q^{M'} y + T_i \in \omega^{M'} (S).\]
Thus there is a coincidence after the $M'$-iteration of the overlap $(R, y, S)$.
In particular, since $\Xi(\Lb)$ is a Meyer set, there are finite equivalence classes of overlaps and so
there exists $l \in \Z_+$ such that for each overlap $\mathcal{O}$ in $\Tk$, $Q^l (\mathcal{O})$ contains a coincidence. Therefore $\Tk$ admits an overlap coincidence.
\qed

\begin{remark} \label{legality-repetitivity}
Note that the legality of every $\Lb$-cluster in a primitive substitution Delone multi-colour set $\Lb$ implies the repetitivity of the associated substitution tiling $\Lb + \Ak$ of $\Lb$ and vice versa.
\end{remark}

Combining the results of Cor\,\ref{purePoint-OverlapCoincidence}, Prop\,\ref{overlap-to-algebraic} and
Prop\,\ref{algebraic-overlap-coincidence}, we get the following theorem.

\begin{theorem} \label{alg-over-coincidence}
Let $\Lb$ be a primitive substitution Delone multi-colour set such that every $\Lb$-cluster is legal and $\Lb$ has FLC.
Then the following are equivalent:
\begin{itemize}
\item[(1)] $(X_{\Lbs}, \R^d, \mu)$ has a pure point dynamical spectrum;
\item[(2)] $\Lb$ admits an algebraic coincidence.
\end{itemize} \qed
\end{theorem}

\section{Algebraic Coincidence to Inter Model Sets} \label{construction-CPS}

\noindent
\subsection{The $Q$-topology} \label{the-Q-topology}
Let $\Lb$ be a primitive substitution Delone multi-colour set.
Define
$$L := \langle \Lam_j \rangle_{j \le m}$$ the group generated by $\Lam_j$, $j \le m$, and let
$$\mathcal{K} := \{x \in \R^d : \Lb + x = \Lb\}$$
be the set of periods of $\Lb$.
Under the assumption that $\Lb$ admits an algebraic coincidence, we introduce a topology on $L$ and find a completion $H$ of the topological group $L$ such that the image of $L$ is a dense subgroup of
$H$. This enables us to construct a cut and project scheme (CPS) such that each point set $\Lam_i$, $i \le m$, arises from the CPS.
In the following lemma we show that the system $\{\alpha + Q^n \Xi(\Lb) + \mathcal{K}: n \in \Z_+, \alpha \in L \}$ satisfies the topological properties for the group $L$ to be a topological group (\cite{HN}, \cite{Bour} and \cite{HR}).

\begin{lemma} \label{Q-topological group}
Let $\Lb$ be a primitive substitution Delone multi-colour set with
expansive map $Q$ such that every $\Lb$-cluster is legal. Suppose that $\Lb$ admits an algebraic coincidence.
Then the system $\{\alpha + Q^n \Xi(\Lb) + \mathcal{K}: n \in \Z_+, \alpha \in L \}$ serves as a neighbourhood base of the topology on $L$ relative to which $L$ becomes a topological group.
\end{lemma}

\noindent{\sc Proof.} From the assumption of an algebraic
coincidence there exist $M \in \Z_+$ and $\xi \in \Lam_i$ such
that \be \label{algebraic-coincidence-formula} Q^M \Xi(\Lb)
\subset \Lam_i - \xi \ \ \ \mbox{for some}\ i \le m. \ee

Consider the family $\mathcal{U} = \{Q^n \Xi(\Lb) + \mathcal{K} \subset L : n \in \Z_+\}$.
We first note that every finite subfamily of $\mathcal{U}$ has a nonempty intersection.
Next we will show that $\mathcal{U}$ satisfies the following property : for every $U \in \mathcal{U}$ and $x \in U$, there exist $V \in \mathcal{U}$ and $V' \in \mathcal{U}$ such that
$V + V \subset U$ and $x + V' \subset U$.
Other properties for a prebase of neighbourhoods of the identity are rather trivial in the Abelian group $L$.
First, note that
\be \label{top-group-con}
Q^M \Xi(\Lb) + Q^M \Xi(\Lb) & \subset & \Xi(\Lb).
\ee
Choose arbitrary $Q^n \Xi(\Lb) + \mathcal{K} \in \mathcal{U}$.
Then with $V:= Q^{M+n}\Xi(\Lb) + \mathcal{K}$
\be
V + V = Q^{M+n} \Xi(\Lb) +\mathcal{K} + Q^{M+n} \Xi(\Lb) +\mathcal{K} &\subset& Q^n \Xi(\Lb) +\mathcal{K}.
\ee
Second,
let $x = Q^n(\alpha - \beta) + k \in Q^n \Xi(\Lb) + \mathcal{K}$, where $\alpha, \beta \in \Lam_j$ for some $j \le m$ and $k \in \mathcal{K}$.
Since every $\Lb$-cluster is legal, there exist $k \in \Z_+$ and $a \in \Lam_j - \Lam_j$ such that
the cluster $\{\alpha, \beta\}$ satisfies $a + \{\alpha, \beta\} \subset \Phi^k(\xi) \cap \Lam_j$ for $\xi \in \Lam_i$ as in
(\ref{algebraic-coincidence-formula}).
So we can find $g \in (\Phi^k)_{ji}$ such that $g(\xi) = a + \alpha$.
From (\ref{algebraic-coincidence-formula}), $g(\xi + Q^M \Xi(\Lb)) \subset \Lam_j$.
So $a + \alpha + Q^{M+k} \Xi(\Lb) \subset \Lam_j$.
Since $a + \beta \in \Lam_j$,
\[ \alpha - \beta + Q^{M+k}\Xi(\Lb) = a + \alpha - (a + \beta) + Q^{M + k} \Xi(\Lb) \subset \Lam_j - (a + \beta) \subset \Lam_j - \Lam_j.
\]
Therefore with $V' := Q^{M+k+n} \Xi(\Lb) + \mathcal{K}$
\[ x + V' = x + Q^{M+k+n} \Xi(\Lb) + \mathcal{K} = Q^n(\alpha - \beta) + Q^{M+k+n} \Xi(\Lb) + \mathcal{K} \subset Q^n \Xi(\Lb) + \mathcal{K}.
\]
Therefore the system $\{\alpha + Q^n \Xi(\Lb) + \mathcal{K} : n \in \Z_+, \alpha \in L \}$ serves as a prebase of neighbourhoods of the topology on $L$ relative to which $L$ becomes a topological group.
In fact the system becomes a neighbourhood base for the topology, since for any $n', n \in \Z_+$ with $n' \ge n$,
$$(Q^{n'} \Xi(\Lb) + \mathcal{K}) \, \cap \, (Q^{n} \Xi(\Lb) + \mathcal{K}) = Q^{n'} \Xi(\Lb) + \mathcal{K} \, \in \, \mathcal{U}.$$
\qed

\medskip

We call the topology on $L$ with the neighbourhood base $\{\alpha
+ Q^n \Xi(\Lb) + \mathcal{K} : n \in \Z_+, \alpha \in L \}$ {\em
$Q$-topology}.


\subsection{Construction of a CPS}

Let $L' = L/\mathcal{K}$ where $L$ and $\mathcal{K}$ are defined
as in Subsec.\,\ref{the-Q-topology}. From \cite[III. \S 3.4, \S
3.5]{Bour} and Lemma 4.1, we know that there exists a complete
Hausdorff topological group of $L'$, which we denote by $H$, for
which $L'$ is isomorphic to a dense subgroup of the complete group
$H$ (see \cite{BM} and \cite{LM2}). Furthermore there is a
uniformly continuous mapping $\psi : L \to H$ which is the
composition of the canonical injection of $L'$ into $H$ and the
canonical homomorphism of $L$ onto $L'$ for which $\psi(L)$ is
dense in $H$ and the mapping $\psi$ from $L$ onto $\psi(L)$ is an
open map, the latter with the induced topology of the completion
$H$. One can directly consider $H$ as the Hausdorff completion of
$L$ vanishing $\mathcal{K}$.

\begin{theorem} \label{model-set-construction}
Let $\Lb$ be a primitive substitution Delone multi-colour set with expansive map $Q$.
Suppose that $\Lb$ admits an algebraic coincidence.
Then there exists a CPS with the locally compact Abelian group $H$ for an internal space
such that for each $j \le m$,
$\Lam_j = \Lam(V_j)$ where $\overline{V_j}$ is a compact set in $H$.
\end{theorem}

\noindent {\sc Proof.} We claim that for any $n \in \Z_{\ge 0}$,
$Q^n \Xi(\Lb) + \mathcal{K}$ is precompact. The argument is
familiar from \cite{BM} and \cite{LM2}, but we provide the
argument for the completeness. Note that $\{Q^{n'} \Xi(\Lb)
+\mathcal{K} : n' \in \Z_{\ge 0}\}$ is a neighbourhood basis for
$0$. We will show that for any $n' \in \Z_{\ge 0}$, $Q^n \Xi(\Lb)
+ \mathcal{K}$ can be covered by some finite translations of
$Q^{n'} \Xi(\Lb) + \mathcal{K}$. For any $n' \in \Z_{\ge 0}$ with
$n' \ge n$, there exists a compact set $C \subset \R^d$ for which
$Q^{n'} \Xi(\Lb) + C = \R^d$, since $Q^{n'}\Xi(\Lb)$ is relatively
dense. So for any $t \in Q^n \Xi(\Lb)$, $t = s + c$ where $s \in
Q^{n'} \Xi(\Lb)$ and $c \in C$. Thus
\[c = t-s \in Q^n \Xi(\Lb) - Q^{n'} \Xi(\Lb) \subset
Q^n \Xi(\Lb) - Q^n \Xi(\Lb).\]
From the assumption of an algebraic coincidence, there exist $M \in \Z_+$ and $\xi \in \Lam_i$ such that $Q^M \Xi(\Lb) \subset \Lam_i - \xi$ for some $i \le m$.
So we get
\be
Q^{n+M} \Xi(\Lb) - Q^{n+M} \Xi(\Lb) &\subset& \Xi(\Lb).
\ee
Thus $Q^M c \in \Xi(\Lb)$.
Since $\Xi(\Lb)$ is discrete, $F:= \Xi(\Lb) \cap Q^M C$ is finite and $Q^M c \in F$.
Thus $t = s+c \in Q^n \Xi(\Lb) + Q^{-M}F$ and we obtain that
\[Q^n \Xi(\Lb) + \mathcal{K} \subset Q^{n'} \Xi(\Lb) + \mathcal{K} + Q^{-M} F.\]
Therefore $Q^n \Xi(\Lb) + \mathcal{K}$ is precompact. This implies
that for each $n \in \Z_{\ge 0}$, \be \label{Xi-is-compact}
\overline{\psi(Q^n \Xi(\Lb) + \mathcal{K})} \ \ \ \mbox{is
compact}. \ee Define $\widetilde{L} := \{(t, \psi(t)) \in \R^d
\times H : t \in L \}$. Applying the same argument as in
\cite[Sec.\,3]{BM}, we note that $\widetilde{L}$ is a uniformly
discrete and relatively dense subgroup in $\R^d \times H$. Then we
can construct a cut and project scheme \be
\begin{array}{ccccc} \label{cut-and-project1}
 \R^{d} & \stackrel{\pi_{1}}{\longleftarrow} & \R^{d} \times H & \stackrel{\pi_{2}}
{\longrightarrow} & H \\
 && \cup \\
L & \longleftarrow & \widetilde{L} & \longrightarrow & \psi(L)\\
 &&  \\
t & \longleftarrow & (t, \psi(t)) & \longrightarrow & \psi(t)
\end{array}
\ee
where $\pi_1$ and $\pi_2$ are the canonical projections.
It is easy to see that ${\pi_1}|_{\widetilde{L}}$ is injective and $\pi_2(\widetilde{L})$ is dense in $H$.
We note that for each $j \le m$, $\Lam_j = \Lam_j + \mathcal{K}$.
So $\Lam_j = \Lam(\psi(\Lam_j + \mathcal{K}))$ and $\overline{\psi(\Lam_j + \mathcal{K})}$ is compact in $H$ from (\ref{Xi-is-compact}).  \qed

\subsection{Existence of Model sets}

In the proof of the following
Lemma\,\ref{find-one-generating-point} and
Prop.\,\ref{genetic-point-set} we make use of the
representability, identifying a primitive substitution Delone
multi-colour set $\Lb$ with the associated substitution tiling
$\Tk = \Lb + \Ak$ as in Subsec.\,\ref{RepresentabilityAsTiling}
where $\Ak$ is the set of tiles arising from the solution of the
adjoint system of equations. Since $\Phi$ is primitive, we can
assume that the substitution matrix $S(\Phi)$ is positive,
replacing $\Phi$ by a power of $\Phi$ if necessary.

In the next Prop.\,\ref{genetic-point-set} we show that if $\Lb$ admits an algebraic coincidence, there exists $\Gb \in X_{\Lbs}$ which is generated from one point and admits an algebraic coincidence at the generating point. This enables us in Prop.\,\ref{genetic-window-open} to show that each point set $\Gam_i$ is a model set whose window is open in $H$ in
CPS (\ref{cut-and-project1}).

\medskip
The following lemma is auxiliary to Prop.\,\ref{genetic-point-set}.

\begin{lemma} \label{find-one-generating-point}
Let $\Lb$ be a primitive substitution Delone multi-colour set with MFS $\Phi$
such that every $\Lb$-cluster is legal.
Suppose that $\supp(\eta' + T_j) \subset (\supp(\omega^N(\eta + T_j)))^{\circ}$ for some $\eta, \eta' \in \Lam_j$ and $N \in \Z_+$, and $f \in (\Phi^N)_{jj}$ for which $f(\eta) = \eta'$.
Then there exists $\Gb = \lim_{n \to \infty} (\Phi^N)^n (y) \in X_{\Lbs}$ for some fixed point $y$ of $f$.
\end{lemma}

\noindent {\sc Proof.}
In this lemma we show how to find a substitution Delone
multi-colour set in $X_{\Lbs}$ which is generated by a fixed
point. We can find a fixed point $y \in \R^d$ of $f$, since $f$ is
an affine map with expansive linear part. Note that
$$ \omega^N(y + T_j) = \{x + T_i : x \in (\Phi^N)_{ij}(y), i \le m\}
$$
(see \cite[Th.\,3.7]{LMS2} for the details). So \be
\label{fixed-tile-in-supertile} y + T_j \in \omega^N(y + T_j)\,.
\ee Notice that $\omega^N(y + T_j)$ is translationally equivalent
to $\omega^N(\eta + T_j)$ and so the relative location of $y +
T_j$ in $\omega^N(y + T_j)$ is same as the relative location of
$\eta' + T_j$ in $\omega^N(\eta + T_j)$, since $f(\eta)=\eta'$ and
$f(y)=y$. So from the assumption that $\supp(\eta' + T_j) \subset
(Q^N(\supp(\eta + T_j)))^{\circ}$, \be
\label{supp-tile-in-interior-supertile} \supp(y + T_j) \subset
(Q^N {\supp(y + T_j)})^{\circ} \,. \ee Now we claim that $0 \in
(\supp(y + T_j))^{\circ}$. By
(\ref{supp-tile-in-interior-supertile}), it is enough to show that
\be \label{zero-in-supp-closure} 0 \in \supp(y + T_j)  \,. \ee In
fact, for any open neighbourhood $U$ of $0$, there exists $s \in
\Z_+$ such that
$$\supp(y + T_j) \subset Q^s U\,.$$
From (\ref{fixed-tile-in-supertile}),
$$Q^{sN} U \cap Q^{sN} \supp(y + T_j) \supset Q^{sN} U \cap \supp(y+T_j)\neq \emptyset \,.$$
Thus
$$U \cap \supp(y + T_j) \neq \emptyset$$ and so the claim is proved.
Let $\Tk' := \lim_{n \to \infty} (\omega^N)^n (y + T_j)$.
Since the generating tile $y + T_j$ contains $0$ in the interior, $\Tk'$ covers $\R^d$ and it is a tiling.
Let
$$\Gb := \lim_{n \to \infty} (\Phi^N)^n (y)\,.$$
Then $\Gb$ is a primitive
substitution Delone multi-colour set in $X_{\Lbs}$ generated from $y \in \Gam_j$.
\qed

\medskip

In the following proposition we show that if the substitution
Delone multi-colour set $\Gb$ obtained from
Lemma\,\ref{find-one-generating-point} admits an algebraic
coincidence at the point $y$, using the assumption of the
algebraic coincidence of $\Lb$.

\begin{prop} \label{genetic-point-set}
Let $\Lb$ be a primitive substitution Delone multi-colour set with expansive map $Q$ and MFS $\Phi$ for which
every $\Lb$-cluster is legal. Suppose that $\Lb$ admits an algebraic coincidence.
Then there exists $\Gb = \lim_{n \to \infty} (\Phi^N)^n (y) \in X_{\Lbs}$ such that
$y + Q^N \Xi(\Lb) \subset \Gam_j$ for some $y \in \Gam_j$, $j \le m$ and $N \in \Z_+$.
\end{prop}

\noindent {\sc Proof.} By the assumption of algebraic coincidence,
there exist $M \in \Z_+$ and $\xi \in \Lam_i$ such that \be
\label{assumption-alg-coincidence} \xi + Q^M \Xi(\Lb) \subset
\Lam_i \ \ \ \mbox{for some} \ i \le m. \ee We again consider the
associated substitution tiling $\Tk := \Lb + \Ak$ of $\Lb$, where
$\Ak =\{T_1, \dots, T_m\}$, and let $A_i = \supp(T_i)$ for $i \le
m$.

It is shown in \cite{lawa} that if $\Lb$ is a substitution Delone
multi-colour set, then there is a finite multiset (cluster) $\Pb
\subset \Lb$ for which $\Lb = \lim_{n \to \infty} \Phi^n(\Pb)$. So
we can find $\eta \in \Lam_j$ for some $j \le m$ and $M' \in \Z_+$
with $M' \ge M$ such that $\eta + T_j $ is fixed under $\omega$
and $\omega^{M'}(\eta + T_j)$ contains $\xi + T_i$. By the
primitivity, we can choose a $j$-type tile $\eta' + T_j $ in the
interior of $\omega^K(\xi + T_i)$ with some $K \in \Z_+$. So
$$\supp(\eta' + T_j ) \subset (\supp(\omega^K(\xi + T_i)))^{\circ} \subset
(\supp(\omega^{M' + K}(\eta + T_j )))^{\circ} \subset (\supp(\omega^{N}(\eta + T_j )))^{\circ},$$
where $N = 2(M' + K)$.
From Lemma \ref{find-one-generating-point},
there exists
\be \label{pointset-generated-from-onepoint}
\Gb = \lim_{n \to \infty} (\Phi^N)^n (y) \in X_{\Lbs}
\ee
for some fixed point $y$ of $f$ where $f(\eta)=\eta'$ and $f \in (\Phi^N)_{jj}$.
Let $\Tk' = \Gb + \Ak$.

We are going to show that $\Gb$ admits an algebraic coincidence at
$y$ using the algebraic coincidence at $\xi$ for $\Lb$ and using
the repetitivity of $\Tk$. Note that $\eta' \in \Phi^K(\xi)$,
which means that there exists $h \in (\Phi^K)_{ji}$ such that
$\eta' = h(\xi)$. Then applying MFS $\Phi$ to the inclusion
(\ref{assumption-alg-coincidence}), we get $h(\xi + Q^M \Xi(\Lb))
\subset h(\Lam_i) \subset \Lam_j$. Thus $\eta' + Q^{M+K} \Xi(\Lb)
\subset \Lam_j$. Then(\ref{assumption-alg-coincidence}) we get
$$\eta' + Q^{M+K} \Xi(\Lb) \subset \Lam_j.$$
Since $Q \Xi(\Lb) \subset \Xi(\Lb)$, $\eta' + Q^{M'+K} \Xi(\Lb) \subset \Lam_j$.
Thus
\[Q^{M' + K} \Xi(\Lb) - Q^{M' + K} \Xi(\Lb) \subset \Xi(\Lb).\]
Let $N = 2(M' + K)$. So
\be \label{Algebraic-condition2}
\eta' + Q^{N} \Xi(\Lb) - Q^{N} \Xi(\Lb) \subset \Lam_j.
\ee
Note that $\Xi(\Tk')$ is also a Meyer set. So as in (\ref{finite-patch})
we can find $a_1, \dots, a_S \in \Xi(\Tk')$ such that for any $a \in \Xi(\Tk')$,
\be \label{Meyer-finite-patches}
 \Tk' \cap ( y + A_j  -a) + a = \Tk' \cap ( y + A_j -a_s) + a_s  \ \ \ \mbox{for some} \ s \le S,
\ee
where $y \in \Gam_j$ as in (\ref{pointset-generated-from-onepoint}).
There exists $p \in \Z_+$ that $(\omega^{N})^p (y + T_j )$ contains the patches of
$\Tk' \cap (y + A_j - a_1), \dots, \Tk' \cap (y + A_j - a_S)$ from (\ref{pointset-generated-from-onepoint}).
Since $\Tk$ is repetitive, there is $r \in \R^d$ such that
\be \label{tile-patch-contained-in-tiling}
\Tk' \cap (y + A_j  - a_s) + r \subset (\omega^{N})^p (y + T_j ) + r \subset \Tk \ \ \ \mbox{for all} \ s \le S.
\ee
Note that $y + r + T_j, \eta + T_j \in \Tk$. Take any $s \le S$.
Since $y + r - \eta \in \Xi(\Tk)$, by (\ref{Algebraic-condition2}) we obtain
\be \label{tile-inclusion-I}
\eta' + T_j  + Q^{N}(y + r - \eta) - Q^{N} a_s \in \Tk \,.
\ee
Now we want to show that
the tile on the left hand side of (\ref{tile-inclusion-I}) is in a patch
$\omega^{N}(\Tk' \cap (y + A_j -a_s)+ r)$.
Since $ \eta' + T_j \in \omega^{N}(\eta + T_j )$,
$$
\eta' + T_j  + Q^{N}(y + r - \eta) - Q^{N} a_s \in \omega^{N}(y + T_j  + r - a_s).
$$
Moreover $$\supp(\omega^N (y + T_j  + r - a_s)) \subset \supp(\omega^N (\Tk' \cap (y + A_j -a_s) + r))$$ and
$$
\omega^N ((\Tk' \cap (y + A_j  -a_s)) + r) \subset \Tk
$$
from (\ref{tile-patch-contained-in-tiling}).
So we can see that
\be
\eta' + T_j + Q^{N}(y + r - \eta) - Q^{N} a_s \in \omega^{N}(\Tk' \cap (y + A_j -a_s)+ r)\,.
\label{containment-of-the-patches}
\ee
Let $f: x \mapsto Q^N x + e$, where $e \in \R^d$. Then we have the identities $\eta' = Q^{N} \eta + e$ and $y = Q^{N} y + e$.
Applying these identities to (\ref{containment-of-the-patches}), we get, for all $s \le S$,
$$y + T_j \in \omega^{N}(\Tk' \cap (y + A_j  -a_s) + a_s)\,.$$
Hence through (\ref{Meyer-finite-patches}), for any arbitrary $a \in \Xi(\Lb)$,
$$y + T_j \in \omega^{N}(\Tk' \cap (y + A_j -a) + a) \,.$$
Thus $y + T_j  - Q^{N}a \in  \omega^N(\Tk' \cap (y + A_j - a)) \subset \Tk'$.
Since $a$ is arbitrary in $\Xi(\Gb)$, $y + Q^{N}\Xi(\Gb) \subset \Gam_j$. \qed

\medskip

\begin{prop} \label{genetic-window-open}
Let $\Lb$ be a primitive substitution Delone multi-colour set.
If $\Lb = \lim_{n \to \infty} \Phi^n (y)$ where $y + Q^M \Xi(\Lb) \subset \Lam_j$ and
$y \in \Lam_j$ for some $M \in \Z_+$ and $j \le m$.
Then for each $i \le m$, $\Lam_i = \Lam(U_i)$ in CPS (\ref{cut-and-project1}) where $U_i$ is an open set and $\overline{U_i}$ is compact in the internal space $H$, i.e. $\Lam_i$ is a model set with an open window.
\end{prop}

\noindent{\sc Proof.} For each $i \le m$ and $z \in \Lam_i$, there
exists $n \in \Z_+$ such that
$$z = Q^n y + e \ \ \  \mbox{for some} \ e \in \R^d,$$
where $f: x \mapsto Q^n x + e$ and $f \in (\Phi^n)_{ij}$.
From $y + Q^M \Xi(\Lb) \subset \Lam_j$, $z + Q^{n + M} \Xi(\Lb) \subset \Lam_i$. Moreover,
$z + Q^{n + M} \Xi(\Lb) + \mathcal{K} \subset \Lam_i$, since $\Lam_i = \Lam_i + \mathcal{K}$.
Thus
$$\Lam_i = \bigcup_{z \in \Lam_i} (z + Q^{M_{z}} \Xi(\Lb) + \mathcal{K}),$$
where $M_{z}$ depends on $z$. Since $\psi$ is an open map from $L$ onto $\psi(L)$, where the latter is with the induced topology of the completion $H$, for each
$z + Q^{M_{z}} \Xi(\Lb) + \mathcal{K}$ there exists an open set $U_{z}$ in $H$ such that
$$\psi(z + Q^{M_{z}} \Xi(\Lb) + \mathcal{K}) = \psi(L) \cap U_{z}.
$$
Since $Ker(\psi) = K$ and $\Lam_i = \Lam_i +K $, $\Lam_i =
\psi^{-1}(\psi(L) \cap U_{i})= \Lam(U_i) $ where $U_i = \bigcup_{z
\in \Lam_i} U_{z}$. Furthermore $\overline{\psi(\Lam_i)} =
\overline{U_i}$ by the denseness of $\psi(L)$ in $H$ and
$\overline{U_i}$ is compact by (\ref{Xi-is-compact}).\qed

\subsection{Two equivalent topologies on $L$} \label{two-topology-on-L}

In this subsection we introduce another topology on $L$ which becomes equivalent to $Q$-topology under the assumption of algebraic coincidence.
Th.\,\ref{geom-gives-modelset} shows a sufficient condition to get inter model set in general setting but the CPS in the theorem is constructed based upon on the new topology. The equivalence of the two topologies gives us the equivalence of the two CPSs. We make use of Th.\,\ref{geom-gives-modelset} to get connection to inter model multi-colour sets for substitution Delone multi-colour set.

Let $\{F_n\}_{n \in \Z_+}$ be a van Hove sequence and let $\Lb',\Lb''$ be two Delone $m$-multi-colour sets
in $\R^d$.
We define
\begin{equation} \label{pseudoMetric}
\rho(\Lb', \Lb'') := \lim_{n \rightarrow \infty} \sup \frac{\sum_{i=1}^m \sharp((\Lam_i' \, \triangle \, \Lam_i'') \cap F_n)}{\Vol(F_n)}.
\end{equation}
Here $\triangle$ is the symmetric difference operator.
Let $P_{\epsilon} = \{x \in L : \rho(x + \Lb, \Lb) < \epsilon \}$ for each $\epsilon > 0$.
From Th.\,\ref{alg-over-coincidence} and \cite[Lemma A.9]{LMS2}, if $\Lb$ admits an algebraic coincidence, then, for any
$\epsilon > 0$, $P_{\epsilon}$ is relatively dense.
In this case the system $\{\alpha + P_{\epsilon} : \epsilon > 0 , \alpha \in L\}$ serves as a neighbourhood base of the topology on $L$ relative to which $L$ becomes a topological group.
We name $P_{\eps}$-topology for this topology on $L$ and denote the space $L$ with $P_{\eps}$-topology by $L_{P}$(see \cite{BM} for
$P_{\eps}$-topology under the name of autocorrelation topology).

\medskip

Let $L_Q$ be the space $L$ with $Q$-topology.
In the following two propositions we show that $L_{Q}$ is topologically isomorphic to $L_{P}$.

\begin{prop} \label{LQ-to-LP}
Let $\Lb$ be a primitive substitution Delone multi-colour set such that every $\Lb$-cluster is legal. Suppose that $\Lb$ admits an algebraic coincidence, then the mapping $\iota : x \mapsto x$ from $L_{Q}$ onto $L_{P}$ is uniformly continuous.
\end{prop}

\noindent{\sc Proof.}
It is enough to show that for each $\epsilon > 0$, there exists $n \in \Z_+$ such that
$Q^n \Xi(\Lb) + \mathcal{K} \subset P_{\epsilon}$.
Let $\Tk = \Lb + \Ak$ be the associated substitution tiling of $\Lb$.
The assumption of algebraic coincidence gives overlap coincidence and, from \cite[Lemma A.9]{LMS2}, there exist $r \in (0,1)$ and $C > 0$ such that for any
$x \in \Xi(\Tk)$
\[ 1 - \mbox{dens}(\Tk \cap (Q^n x + \Tk)) \le C r^n . \]
Since
\[\mbox{dens}(\Tk \cap (Q^n x + \Tk)) = \sum_{i = 1}^m \mbox{freq}(\{T_i, Q^n x + T_i\}, \Tk) \cdot
\mbox{Vol}(A_i)\]
and
\[\mbox{freq}(\{T_i, Q^n x + T_i\}, \Tk) = \mbox{dens}(\Lam_i \cap (Q^n x + \Lam_i)),\]
we get
\begin{eqnarray*}
\lefteqn {1 - \mbox{dens}(\Tk \cap (Q^n x + \Tk))} \\
& = & 1 - \sum_{i =1}^m \mbox{dens}(\Lam_i \cap (Q^n x + \Lam_i)) \cdot \mbox{Vol}(A_i) \\
& = & \sum_{i =1}^m \mbox{dens}(\Lam_i) \cdot \mbox{Vol}(A_i) -
\sum_{i =1}^m \mbox{dens}(\Lam_i \cap (Q^n x + \Lam_i)) \cdot \mbox{Vol}(A_i) \\
& = & \sum_{i =1}^m \frac{1}{2} (\mbox{dens}(\Lam_i \, \triangle \, (Q^n x + \Lam_i))) \cdot
\mbox{Vol}(A_i).
\end{eqnarray*}
Let $V_0 = \mbox{min} \{\mbox{Vol}(A_i): i \le m\}$. Then
\[V_0 \cdot \sum_{i =1}^m \frac{1}{2} (\mbox{dens}(\Lam_i \, \triangle \, (Q^n x + \Lam_i))) \le C r^n .\]
So for all $x \in \Xi(\Lb)$,
\[ \mbox{dens}(\Lb \, \triangle \, (Q^n x + \Lb)) \le C' r^n  \ \ \ \mbox{where} \ C' = \frac{2C}{V_0} > 0.
\]
Then $Q^n x \in P_{C'r^n}$ and thus for any $\epsilon > 0$, we can find $n \in \Z_+$ satisfying $C' r^n  < \epsilon$ so that
$$Q^n \Xi(\Lb) + \mathcal{K} \subset P_{\epsilon}.$$ \qed

\begin{prop} \label{LP-to-LQ}
Let $\Lb$ be a primitive substitution Delone multi-colour set such that every $\Lb$-cluster is legal. Suppose that $\Lb$ admits an algebraic coincidence, then the mapping $\iota^{-1} : x \mapsto x$ from $L_{P}$ onto
$L_{Q}$ is uniformly continuous.
\end{prop}

\noindent{\sc Proof.}
It is enough to show that for any $n \in \Z_+$ there exists $\epsilon > 0$ such that $P_{\epsilon} \subset Q^n \Xi(\Lb) + \mathcal{K}$.
For the associated substitution tiling $\Tk = \Lb + \mathcal{A}$ of $\Lb$.
Since $\Xi(\Tk)$ is a Meyer set from Lemma\,\ref{alg-coincidence-Meyer}, we can find small $\zeta > 0$ such that
$\Tk \cap (r + t +  \Tk) = \emptyset$ for all $t \in \Xi(\Tk)$ and
$r \in B_{\zeta}(0) \, \backslash \, \{0\}$.
From \cite[Lemma\,3.5 and Th.\,2.14]{sol3}, for any $0 < \eps \le \frac{\zeta}{||Q^n||}$ there exists $\delta_{\eps} > 0$ such that
if $d(\Tk-x, \Tk -(x+t)) < \delta_{\eps}$ for some $x \in \R^d$ and $t \in \Xi(\Tk)$, then there exist $g_j, g'_j \in \mathcal{K}$ for $1 \le j \le n$, such that
\[ d(\Tk - Q^{-n}x - \sum_{j=1}^n Q^{-n+j-1}g_j, \Tk - Q^{-n}(x+t) - \sum_{j=1}^n Q^{-n+j-1}g'_j) < \eps\,.
\]
By the metric $d$ on $X_{\Tk}$ of (\ref{metric-tilings}), there exists
$h \in B_{\eps}(0)$ such that
$$\Tk - Q^{-n}x - \sum_{j=1}^n Q^{-n+j-1}g_j - h \ \  \mbox{agrees with} \
\Tk - Q^{-n}(x+t) - \sum_{j=1}^n Q^{-n+j-1}g'_j \ \ \mbox{on} \ B_{1/\eps}(0).
$$
So
$$ \Tk - x - \sum_{j=1}^n Q^{j-1}g_j - Q^n h \ \  \mbox{agrees with} \
\Tk - (x+t) - \sum_{j=1}^n Q^{j-1}g'_j \ \ \mbox{on} \ Q^n B_{1/\eps}(0).
$$
Note that $$t - \sum_{j=1}^n Q^{j-1}g_j + \sum_{j=1}^n Q^{j-1}g'_j \in \Xi(\Tk),$$
since $Q \mathcal{K} \subset \mathcal{K}$ and $\Xi(\Tk) + \mathcal{K} = \Xi(\Tk)$.
Since $||Q^n|| \eps \le \zeta$ and the choice of $\xi$, $Q^n h = 0$. Thus $h = 0$.
So
\[ Q^{-n}t - \sum_{j=1}^n Q^{-n+j-1}g_j + \sum_{j=1}^n Q^{-n+j-1}g'_j \in \Xi(\Tk).\]
Thus
$$ t = Q^n z + w, \ \ \ \mbox{where} \ z \in \Xi(\Tk) \ \mbox{and} \
w = \sum_{j=1}^n Q^{j-1}(g_j - g'_j) \in \mathcal{K},
$$
and hence $t \in Q^n \Xi(\Tk) + \mathcal{K}$.

If $t \in P_{\eps}$, $\rho(t+\Lb,\Lb) < \eps$. This means that for small $\eps > 0$ there is a big area of overlaps in $\R^d$ between $t+\Tk$ and $\Tk$ so that $d(\Tk + x, \Tk + x -t)$ is small for some $x \in \R^d$.
So we can choose small $\epsilon > 0$ so that for any $t \in P_{\epsilon}$, $d(\Tk + x, \Tk + x -t) < \delta_{\epsilon}$
for some $x \in \R^d$ by the definition of $P_{\epsilon}$. Then $t \in Q^n \Xi(\Tk) + \mathcal{K}$. Hence $P_{\epsilon} \subset
Q^n \Xi(\Lb) + \mathcal{K}$.
\qed

\medskip

\begin{remark}
{\em From Prop.\,\ref{LQ-to-LP} and Prop.\,\ref{LP-to-LQ}, $L_{P}$ is topologically isomorphic to $L_{Q}$. Thus the completion of
$L_{P}$ is topologically isomorphic to the completion $H$ of $L_{Q}$. We will identify the former with $H$.
Thus $\phi := \psi \cdot \iota^{-1}  : L_{P} \to H$ is uniformly continuous,
$\phi(L_{P})$ is dense in $H$, and the mapping $\phi$ from $L_{P}$ onto $\phi(L_{P})$ is an open map, the latter with the induced topology of the completion $H$.
Therefore we can consider the CPS (\ref{cut-and-project1}) with an internal space $H$ which is a completion of $L_{P}$.
Note that since $\Lb$ is repetitive, $\bigcap_{\epsilon > 0} P_{\epsilon} = \mathcal{K}$ and
$\mathcal{K} = \overline{\{0\}}$ in $L_{Q}$. }
\end{remark}

\subsection{Inter Model Sets} \label{algebraic-coincidence-modelsets}

\subsubsection{A continuous map between two dynamical hulls}

In this subsection we show that for a primitive substitution Delone multi-colour set $\Lb$ with pure point spectrum, there exists a continuous map from $X_{\Lbs}$ to $\mathbb{A}(\Lb)$. This continuous map was first introduced in \cite{BLM}.


Let us define an {\em autocorrelation group} $\mathbb{A}(\Lb)$.
Let $\widetilde{\mathcal{D}}$ be the set of all Delone $m$-multi-colour sets in $\R^d$. Define an equivalence relation on $\widetilde{\mathcal{D}}$ by
$\Lb' \equiv \Lb'' \Leftrightarrow \rho(\Lb', \Lb'')=0$.
Let $\mathcal{D} := \widetilde{\mathcal{D}}/ \equiv$ and
let $\rho$ also denote the resulting $\R^d$-invariant metric on $\mathcal{D}$.
Now we define a new uniformity on $\mathcal{D}$, which mixes
the autocorrelation topology with the standard topology of $\R^d$
using the sets
$$U(V,\epsilon) = \{(\Lb', \Lb'') \in \mathcal{D} \times \mathcal{D} : \rho(-v + \Lb', \Lb'') < \epsilon \; {\textrm{for some}} \; v \in V \} \, $$
where $\eps > 0$ and $V$ is a neighbourhood of $0$.
Then $\mathcal{D}$ is a complete space \cite{BLM}
and its elements can be identified as Delone multi-colour sets in $\R^d$ up to density $0$ changes.
Notice that the topology induced by this uniformity is not same with the topology ($P_{\eps}$-topology) induced by the
metric $\rho$.
We define $\mathbb{A}(\Lb)$ as the closure of the orbit $\R^d + \Lb$ with the new uniformity in $\mathcal{D}$ (see \cite{BLM}and \cite{LM2} for more about $\mathbb{A}(\Lb)$).

The following theorem is proved in \cite{BLM} in Delone sets with one colour.
The argument can be extended into Delone multi-colour sets without difficulty.

\begin{theorem} \label{continuous-map-beta}
Let $\Lb$ be a Delone multi-colour set in $\R^d$ with UCF such that $\Xi(\Lb)$ is a Meyer set. If the dynamical system $(X_{\Lbs}, \mu, \R^d)$ has pure point spectrum with continuous eigenfunctions, then there exists a continuous $\R^d$-map
$$\beta : X_{\Lbs} \rightarrow \A(\Lb) \,, $$
in which $\beta : \Gb \mapsto \Gb \ \mbox{mod} \ \equiv$.
\qed
\end{theorem}

In substitution Delone multi-colour sets the condition of continuous eigenfunctions is already implicit:

\begin{theorem} \cite[Th.\,2.13]{sol3} \label{conti-eigenfunctions}
If $\Lb$ is a primitive substitution Delone multi-colour set with FLC such that every $\Lb$-cluster is legal.
Then every measurable eigenfunction for the system $(X_{\Lbs}, \mu, \R^d)$ coincides with a continuous function $\mu$-a.e.
\end{theorem}

\medskip

Since a primitive substitution Delone multi-colour set with FLC has UCF (see \cite{LMS2}), we combine Th.\,\ref{continuous-map-beta}, Th.\,\ref{conti-eigenfunctions} and Th.\,\ref{purePoint-Meyer} and get the following corollary.

\begin{cor} \label{sub-ppd-conti}
Let $\Lb$ be a primitive substitution Delone multi-colour set with FLC such that every $\Lb$-cluster is legal. If the dynamical system $(X_{\Lbs}, \mu, \R^d)$ has pure point spectrum. Then there exists a continuous $\R^d$-map
$$\beta : X_{\Lbs} \rightarrow \A(\Lb),$$
in which $\beta : \Gb \mapsto \Gb \ \mbox{mod} \ \equiv$.
\qed
\end{cor}

\subsubsection{Algebraic coincidence to inter model sets}

In this subsection we show that if a substitution Delone multi-colour set admits an algebraic coincidence then it is an inter model multi-colour set.

A continuous $\R^d$-map $\beta : X_{\Lbs} \rightarrow \mathbb{A}(\Lb)$ is called a torus parametrization on $X_{\Lbs}$.\footnote{The terminology of a torus parametrization arises from the model set cases first studied where
(in the set-up that we have here) $\mathbb{A}(\Lb)$ would have been a torus.}  An element
$\Gb \in X_{\Lbs}$ is {\em non-singular} for this parametrization if $\beta^{-1}(\beta(\{\Gb\}))= \{\Gb\}$. The set of non-singular elements of $X_{\Lbs}$ is invariant under the translation action of
$\R^d$.

The result of the following theorem is based on a CPS, taking the completion of $L_{P}$ as an internal space (see \cite{LM2}). Since we have shown that the completion of $L_{P}$ is topologically isomorphic to the completion of $L_{Q}$, we can use the
CPS (\ref{cut-and-project1}) in Th.\,\ref{geom-gives-modelset}.

Notice that $\mathbb{A}(\Lb)$ is isomorphic to a torus $(\R^d \times H)/\widetilde{L}$ by
\cite[Prop.\,3.2]{LM2}.

\begin{theorem} \label{geom-gives-modelset}\cite[Prop.\,4.6]{LM2}
Let $\Lb$ be a multi-colour set in $\R^d$ with repetitivity.
Suppose that there exists a continuous $\R^d$-map $\beta : X_{\Lbs} \rightarrow \mathbb{A}(\Lb)$ and $\Lam({V_i}^{\circ}) \subset \Lam_i \subset \Lam(\overline{V_i})$ where $\overline{V_i}$ is compact, ${V_i}^{\circ} \neq \emptyset$, and $\partial{V_i}$ has empty interior for each $i \le m$ with respect to CPS (\ref{cut-and-project1}).
Then there exists a non-singular element $\Lb'$ in $X_{\Lbs}$ such that
$\Lam'_i = \Lam(W_i)$ where $W_i$ is compact and $W_i = \overline{{W_i}^{\circ}}$ for each $i \leq m$ with respect to the same CPS, and so
for each $\Gb \in X_{\Lbs}$ there exists $(-s, -h) \in \R^d \times H$ so that
\[-s + \Lam(h + {W_i}^{\circ}) \subset \Gam_i \subset -s + \Lam(h + W_i) \ \ \ \mbox{for each} \ i \le m.
\]
In other words, every $\Gb \in X_{\Lbs}$ is an inter model multi-colour set.
\end{theorem}

\medskip

\begin{theorem} \label{alg-coincidence-model-set}
Let $\Lb$ be a primitive substitution Delone multi-colour set such that every $\Lb$-cluster is legal. Suppose that $\Lb$ admits an algebraic coincidence.
Then for each $\Gb \in X_{\Lbs}$ there exists $(-s, -h) \in \R^d \times H$ satisfying
\[-s + \Lam(h + {W_i}^{\circ}) \subset \Gam_i \subset -s + \Lam(h + W_i) \ \ \ \mbox{for each}~ i \le m,\]
where $W_i$ is compact and $W_i = \overline{{W_i}^{\circ}} \neq \emptyset$, with respect to the
CPS (\ref{cut-and-project1}). In other words, every $\Gb \in X_{\Lbs}$ is an inter model multi-colour set.
In particular, $\Lb$ is an inter model multi-colour set.
\end{theorem}

\noindent{\sc Proof.} From Prop.\,\ref{genetic-point-set} and
Prop.\,\ref{genetic-window-open}, there exists $\Lb' \in X_{\Lbs}$
for which each $\Lam'_i = \Lam(V_i)$, where $V_i \neq 0$ is an
open set and $\overline{V_i}$ is compact in $H$. Here we note that
the boundary $\partial{V_i}$ of the open set $V_i$ has empty
interior. Note that algebraic coincidence implies pure point
dynamical spectrum by Th.\,\ref{alg-over-coincidence}. Applying
Cor.\,\ref{sub-ppd-conti} and Th.\,\ref{geom-gives-modelset} to
$\Lb'$, we complete the theorem. \qed

\medskip

\section{Inter model sets to algebraic coincidence } \label{InterModelsets-to-algCoincidence}

\noindent
We will show that if a substitution Delone multi-colour set $\Lb$ is an inter model multi-colour set then $\Lb$ admits an algebraic coincidence.

For each compact set $K \subset \R^d$, we define
$$T_K(\Lb) := \{ t \in L : t + (\Lb \cap K) = \Lb \cap (t + K)\}.$$

\medskip
We prove the following auxiliary lemma for Th.\,\ref{model-set-alg-coincidence}.

\begin{lemma} \label{intermodelset-locatorset}
Let $\Lb$ be a Delone multi-colour set in $\R^d$. If for each $i
\le m$, $\Lam({W_i}^{\circ}) \subset \Lam_i \subset \Lam(W_i)$ for
some compact set $W_i \neq \emptyset$ with $W_i =
\overline{{W_i}^{\circ}}$, in some CPS. Then $T_F(\Lb) - T_F(\Lb)
\subset - \xi + \Lam_j$ for some compact set $F$ and $\xi \in
\Lam({W_i}^{\circ})$.
\end{lemma}

\noindent
{\sc Proof.}
For any $t \in T_K(\Lb)$,
$t + (\Lb \cap K) \subset \Lb \cap (t + K)$ and $t + (\Lb \cap K) \supset \Lb \cap (t + K)$.
So for each $i \le m$,
\be
\psi(t) + \psi(s) \in W_i && \forall \ s \in \Lam_i \cap K  \ \ \ \mbox{and}\nonumber\\
\psi(t) + \psi(s) \notin {W_i}^{\circ} && \forall \ s \in (L \setminus \Lam_i) \cap K. \nonumber
\ee
Let
\[ W_{i, K} := \bigcap \{-\psi(s) + W_i : s \in \Lam_i \cap K\} \setminus \bigcup
\{-\psi(s) + {W_i}^{\circ} : s \in (L \setminus \Lam_i) \cap K\},
\]
for each $i \le m$.
Then we can say that
\[ T_K(\Lb) \subset \Lam(\bigcap_{i \le m} W_{i, K}).\]
Fix any $j \le m$. Since ${W_j}^{\circ} \neq \emptyset$, we can find $\xi$ such that $\xi \in \Lam({W_j}^{\circ}) \subset \Lam_j$.
Since $-\psi(\xi) + {W_j}^{\circ}$ contains a neighbourhood of $0$ and $H$ is a locally compact Abelian group, there is a neighbourhood $U$ of $0$ in $H$ such that $U - U \subset -\psi(\xi) + {W_j}^{\circ}$.
Let
\[I = \{t \in H : t + W_i = W_i~~ \mbox{for all}~i \le m \}.\]
Since $W_i = \overline{{W_i}^{\circ}}$ for all $i \le m$,
\[\{t \in H : t + {W_i}^{\circ} = {W_i}^{\circ}~~ \mbox{for all}~i \le m \} = I.\]
So
\be
(U + I) - (U + I) &=& U -U +I \nonumber \\
 &\subset& -\psi(\xi) + {W_j}^{\circ} + I \nonumber \\
 &=&  -\psi(\xi) + {W_j}^{\circ} \,. \nonumber
\ee
Note that $U + I$ is a neighbourhood of $0$ in $H$.

We claim that
\[\bigcap \{ -\psi(s) + W_i : s \in \Lam_i, i \le m\} = I.\]
Notice first that
$$ 0 \in \bigcap \{ -\psi(s) + W_i : s \in \Lam_i, i \le m\} \neq \emptyset \,.$$
For any $c \in \bigcap \{- \psi(s) + W_i : s \in \Lam_i, i \le m\}$, $\psi(\Lam_i) \subset -c + W_i$ for all $i \le m$. So $\overline{\psi(\Lam_i)} = W_i \subset -c + W_i$ for all $i \le m$.
In fact, $W_i = -c + W_i$ for all $i \le m$, since $W_i - W_i$ is compact
(see \cite[Prop.\,5.2]{LM2} for the detailed proof).
Thus $\bigcap \{- \psi(s) + W_i : s \in \Lam_i, i \le m\} \subset I$.
On the other hand, for any $c' \in I$, $c' + W_i = W_i$ for all $i \le m$, and so
\[-\psi(s) + c' + W_i =
-\psi(s) + W_i ~~\mbox{for all}~ s \in \Lam_i ~\mbox{and}~ i \le m.\]
Since
$0 \in \bigcap \{ -\psi(s) + W_i : s \in \Lam_i, i \le m\}$,
 $$c' \in -\psi(s) + W_i ~~\mbox{for all}~ s \in \Lam_i ~~\mbox{and}~ i \le m.$$
This shows that
$c' \in \bigcap \{ -\psi(s) + W_i : s \in \Lam_i, i \le m\}$.
Therefore the claim is proved.

So now we have $\bigcap \{ (-\psi(s) + W_i) \setminus (U + I) : s
\in \Lam_i, i \le m\} = \emptyset$. Since each $(-\psi(s) + W_i)
\setminus (U + I)$ is compact, by the finite intersection property
for compact sets there is a finite set $F \subset L$ such that $F
\subset \bigcup_{i \le m} \Lam_i$ and
\[\bigcap \{ (-\psi(s) + W_i) \setminus (U + I) : s \in \Lam_i \cap F, \, i \le m\} = \emptyset.
\]
Thus $ U + I \supset \bigcap \{ -\psi(s) + W_i : s \in \Lam_i \cap F, \, i \le m\} \supset
\bigcap_{i \le m} W_{i, F}$.
Then for the compact set $F \subset \R^d$,
\be \label{lambda-has-openset}
T_F(\Lb) - T_F(\Lb) & \subset & \Lam(\bigcap_{i \le m} W_{i,F}) - \Lam(\bigcap_{i \le m} W_{i,F}) \nonumber \\
& \subset & \Lam((U+I)-(U+I)) \nonumber \\
& \subset & \Lam(-\psi(\xi) + {W_j}^{\circ}) \nonumber \\
& \subset & -\xi + \Lam_j.
\ee
\qed

\begin{theorem} \label{model-set-alg-coincidence}
Let $\Lb$ be a primitive substitution Delone multi-colour set such that every $\Lb$-cluster is legal.
Suppose that $\Lb$ is an inter model multi-colour set.
Then $\Lb$ admits an algebraic coincidence.
\end{theorem}

\noindent
{\sc Proof.}
From the assumption there is a following cut and project scheme :
\be \label{cut-and-project2}
\begin{array}{ccccc}
 \R^d & \stackrel{\pi_{1}}{\longleftarrow} & \R^d \times H & \stackrel{\pi_{2}}
{\longrightarrow} & H \,,\\
  && \cup \\
  & & \widetilde{L} & &
\end{array}
\ee
where $H$ is a locally compact Abelian group, $\widetilde{L}$ is a lattice in $\R^d \times H$,
$\pi_1$ and $\pi_2$ are canonical projections,
$\pi_1|_{\widetilde{L}}$ is one-to-one, and $\pi_2(\widetilde{L})$ is dense in $H$.
Let $L = \pi_1(\widetilde{L})$. We define $\psi : L \rightarrow H$ by
$\psi(x) = \pi_2(\pi_1^{-1}(x))$.
Then $s+ \Lam({W_i}^{\circ}) \subset \Lam_i \subset s+ \Lam(W_i)$ for some $s \in \R^d$ and non-empty compact set $W_i$ with
$W_i = \overline{{W_i}^{\circ}}$ for each $i \le m$ with respect to the CPS (\ref{cut-and-project2}).
We can assume that $\Lam({W_i}^{\circ}) \subset \Lam_i \subset \Lam(W_i)$ without loss of generality.

From Lemma\,\ref{intermodelset-locatorset}, we have
\be \label{locatorset-containment}
T_F(\Lb) - T_F(\Lb) \subset - \xi + \Lam_j
\ee
for some compact set $F$ and $\xi \in \Lam({W_i}^{\circ})$.
Let $K = F + \bigcup_{i \le m}(\supp(T_i))$.
Since every $\Lb$-cluster is legal, there exist $\alpha \in \Lam_k$ for some
$k \le m$ and $N \in \Z_+$ satisfying $\Lb \cap K \subset z + \Phi^N(\alpha)$ for some $z \in \R^d$.
Note that for any $\beta \in \Lam_k$,
\[z + Q^N(\beta - \alpha) + \Phi^N(\alpha) = z + \Phi^N(\beta).
\]
Thus $Q^N(\beta - \alpha) + (\Lb \cap K) \subset z + \Phi^N(\beta)$ and so
\[-z + Q^N(\beta - \alpha) + (\Lb \cap K) \subset \Lb.
\]
We note further that
\[-z + Q^N(\beta - \alpha) + (\Lb \cap K) \subset \Lb  \cap (-z + Q^N(\beta - \alpha) + K).
\]
By the choice of $K$ and the fact that $\Lb + \Ak$ is a tiling,
\[-z + Q^N(\beta - \alpha) + (\Lb \cap F) = \Lb  \cap (-z + Q^N(\beta - \alpha) + F).
\]
So $-z + Q^N(\beta - \alpha) \in T_F(\Lb)$. This shows that $-z +
Q^N(\Lam_k - \alpha) \in T_F(\Lb)$. Thus $Q^N(\Lam_k - \Lam_k)
\subset T_F(\Lb) - T_F(\Lb)$. By (\ref{locatorset-containment})
and the fact that $Q^{N+l}\Xi(\Lb) \subset Q^N(\Lam_k - \Lam_k)$
for some $l \in \Z_+$ by the primitivity of the substitution,
there exist $M \in \Z_+$ and $\xi \in \Lam_j$ such that $Q^M
\Xi(\Lb) \subset \Lam_j - \xi$. This completes the proof. \qed

\medskip

The following theorem states the main result of this paper.

\begin{theorem} \label{main-theorem}
Let $\Lb$ be a primitive substitution Delone multi-colour set such that every $\Lb$-cluster is legal and $\Lb$ has FLC.
Then the following are equivalent:
\begin{itemize}
\item[(1)] $\Lb$ has pure point dynamical spectrum; \item[(2)]
$\Lb$ admits an algebraic coincidence; \item[(3)] $\Lb$ is an
inter model multi-colour set.
\end{itemize}
\end{theorem}

\noindent
{\sc Proof.} The proof goes as follows:
\begin{itemize}
\item[]$(1) \Leftrightarrow (2)$ by
Th.\,\ref{alg-over-coincidence}.
\item[]$(2) \Leftrightarrow (3)$
by Th.\,\ref{alg-coincidence-model-set} and
Th.\,\ref{model-set-alg-coincidence}. \qed
\end{itemize}

\medskip

Any tiling $\Tk$ can be converted into a Delone multiset by simply
choosing a point for each tile so that the chosen points for tiles
of the same type are in the same relative position in the tiles.
So we give a corresponding result of Th.\,\ref{main-theorem} on
substitution tilings.

The following lemma is taken from a lecture note of Boris
Solomyak. We provide the proof here, since there is no direct
reference for it. One can see the similar arguments in
\cite{lawa}, \cite{sol4}, and \cite{LMS2}.

\begin{lemma}
Let $\Tk$ be a repetitive fixed point of a primitive substitution
such that $ \Tk = \bigcup_{j=1}^m (T_j + \Lam_j)$. Then $\Lb_{\Tk}
:= (\Lam_i)_{i \le m}$ is a primitive substitution Delone multiset
and every $\Lb_{\Tk}$-cluster is legal.
\end{lemma}

\noindent {\sc Proof.} Let $\omega$ be the corresponding
tile-substitution for $\Tk$. Then \[ \Tk = \bigcup_{j=1}^m
(\omega(T_j) + Q \Lam_j) = \bigcup_{j=1}^m \left(
\bigcup_{i=1}^m(T_i + \mathcal{D}_{ij}) + Q \Lam_j \right) =
\bigcup_{i=1}^m \left( T_i + \bigcup_{j=1}^m(Q \Lam_j +
\mathcal{D}_{ij}) \right)\,. \] Thus
\[\Lam_i = \bigcup_{j=1}^m (Q \Lam_j + \mathcal{D}_{ij}), \ \
 \ i \le m\,.\]
Every $\Lb_{\Tk}$-cluster is legal from
Remark\,\ref{legality-repetitivity}. \qed

\begin{theorem}
Let $\Tk$ be a repetitive fixed point of a primitive substitution
with FLC. Then the following are equivalent:
\begin{itemize}
\item[(1)] $\Tk$ has pure point dynamical spectrum; \item[(2)]
$\Tk$ admits an overlap coincidence; \item[(3)] $\Lb_{\Tk}$ is an
inter model multi-colour set.
\end{itemize}
\end{theorem}

\noindent {\sc Proof.} The proof goes as follows:
\begin{itemize}
\item[] $(1) \Leftrightarrow (2)$ by
Cor.\,\ref{purePoint-OverlapCoincidence}. \item[] $(2)
\Leftrightarrow (3)$ by Prop.\,\ref{overlap-to-algebraic},
Prop.\,\ref{algebraic-overlap-coincidence}, and
Th.\,\ref{main-theorem}. \qed
\end{itemize}

\section{Further study}

\noindent When the legality in Th.\,\ref{main-theorem} is dropped,
finding the corresponding substitution tilings is not obvious.
However with an assumption of repetitivity of $\Lb$, we get a
notion of multi-tilings which is introduced in \cite{Lee}. Can
Th.\,\ref{main-theorem} be extended when $\Lb$ is assumed to be
only repetitive?

In lattice substitution Delone multi-colour sets, modular coincidence was introduced as a condition equivalent to pure point diffractivity, and it was proved to be computable (see \cite{LM1} and \cite{LMS2}).
Is there an algorithm for checking algebraic coincidence in substitution Delone multi-colour set?

There is a considerable amount of ongoing work on Pisot-type
substitution sequences for the study of number theory, discrete
geometry, geometrical combinatorics, mathematical quasicrystals
and spectral theory. As a special case of substitutions there are
Pisot substitutions in $1$-dimension each of whose substitution
matrices has one eigenvalue strictly bigger than $1$ and other
eigenvalues strictly between $0$ and $1$ in modulus. It has been
conjectured that every Pisot substitution dynamical system in
$1$-dimension has pure point spectrum. Here the algebraic
coincidence is an alternative way to determine pure point spectrum
in the Pisot substitutions. Throughout correspondence with Valerie
Berthe, it is noted that the algebraic coincidence is necessary if
an exclusive inner point exists, which is conjectured to hold for
every Pisot unit substitution in 1-dimension (see \cite{Akiyama}).
Does every Pisot unit substitution admit algebraic coincidence?
Bernd Sing's thesis \cite{Sing} deals with this problem and
provides many equivalence properties to the algebraic coincidence.

Although we have shown in this paper the equivalence between the
notions of inter model sets and pure point spectrum in
substitution Delone multi-colour sets, we do not know the measure
of the boundary of the window of the inter model set. When the
underlying structure of substitution Delone multi-colour set is on
a lattice, we know that the measure of the boundary is zero from
\cite{LM1} and \cite{LMS2}. It is a remaining question that in
substitution Delone multi-colour sets (not assumed to be on
lattices) or in Delone multi-colour sets if there is any inter
model set with pure point spectrum whose window has boundary of
non-zero measure.

In \cite{BLM} regular model sets are characterized in terms of
their associated dynamical systems. Can inter model sets be
characterized in terms of their associated dynamical systems as
well? Would it be possible to extend the equivalence of inter
model sets and pure point spectra for general model sets (not
assumed to be substitution Delone sets)?

\section{Acknowledgment}

\noindent
The author is grateful to Valerie Berthe, Robert V. Moody, Boris Solomyak, and Nicolae Strungaru for helpful discussions and insight.
The main ideas on this paper have been established during her Ph.D program. She is indebted to Robert V. Moody for his guidance and Boris Solomyak for valuable comments and encouragements.

\end{document}